\documentclass[a4paper,12pt]{article}

\makeatletter
\@addtoreset{footnote}{page}
\makeatother

\usepackage{amsthm}
\usepackage{amsmath,amssymb,latexsym,amsfonts,mathrsfs}
\usepackage{color}

\renewcommand{\d}{{\rm d}} 

\newcommand{\e}{{\rm e}} 

\newcommand{\law}{\stackrel{{\rm law}}{=}}

\newcommand{\tend}[2]{\mathrel{\mathop{\longrightarrow}\limits^{#1}_{#2}}}

\renewcommand{\tilde}{\widetilde}

\newcommand{\up}{\ensuremath{\uparrow}}
\newcommand{\down}{\ensuremath{\downarrow}}

\newcommand{\rbra}[1]{\!\left( #1 \right)} 
\newcommand{\cbra}[1]{\!\left\{ #1 \right\}} 
\newcommand{\sbra}[1]{\!\left[ #1 \right]} 

\newcommand{\bP}{\ensuremath{\mathbb{P}}}
\newcommand{\bQ}{\ensuremath{\mathbb{Q}}}

\newcommand{\cF}{\ensuremath{\mathcal{F}}}

\newcommand{\cL}{\ensuremath{\mathcal{L}}}

\newcommand{\cP}{\ensuremath{\mathcal{P}}}

\newcommand{\ve}{\ensuremath{\mbox{{\boldmath $e$}}}}

\newcommand{\vn}{\ensuremath{\mbox{{\boldmath $n$}}}}

\theoremstyle{plain}
\newtheorem{Thm}{Theorem}[section]
\newtheorem{Lem}[Thm]{Lemma}

\theoremstyle{definition}

\newtheorem{Rem}[Thm]{Remark}

\newcommand{\Proof}[2][Proof]{\begin{proof}[{#1}] #2 \end{proof}}

\numberwithin{equation}{section}

\makeatletter
\renewcommand\section{\@startsection {section}{1}{\z@}%
                                   {-3.5ex \@plus -1ex \@minus -.2ex}%
                                   {2.3ex \@plus.2ex}%
                                   {\normalfont\large\bf}}
\makeatother

\makeatletter
\renewcommand\subsection{\@startsection {subsection}{1}{\z@}%
                                   {-3.5ex \@plus -1ex \@minus -.2ex}%
                                   {2.3ex \@plus.2ex}%
                                   {\normalfont\normalsize\bf}}
\makeatother

\begin{document}

\begin{center}
{\bf \Large 
Local time penalizations with various clocks for one-dimensional diffusions
}
\end{center}

\begin{center}
Christophe Profeta\footnotemark, 
Kouji Yano\footnotemark\,and Yuko Yano\footnotemark 
\end{center}\footnotetext[1]{Universit\'e d'Evry Val-d'Essonne; \footnotemark[2]Kyoto University; \footnotemark[3]Kyoto Sangyo University. \footnotemark[1]\footnotemark[2]These authors were supported by JSPS-MAEDI Sakura program. \footnotemark[2]This author was supported by MEXT KAKENHI grant no.'s 26800058, 24540390 and 15H03624. \footnotemark[3]This author was supported by MEXT KAKENHI grant no. 23740073.}

\begin{abstract}
We study some limit theorems for the law of a generalized  one-dimensional diffusion weighted and normalized by a non-negative function of the local time  evaluated at a parametrized family of random times (which we will call a {\em clock}).
As the clock tends to infinity, we show that the initial process converges towards a new penalized process, which generally depends on the chosen clock. However, unlike with deterministic clocks, no specific assumptions are needed on the resolvent of the diffusion. We then give a path interpretation of these penalized processes via some universal $ \sigma $-finite measures.

\end{abstract}

\section{Introduction}

\subsection{Penalizations}
The systematic studies of penalizations started in 2003 
with the works of Roynette, Vallois and Yor; see for instance \cite{MR2229621}, \cite{MR2253307} for early papers, and \cite{RoyYor} for a monograph on this subject. 
The starting point of our study is the following classical penalization result: 
if $|B|$ is a standard reflected Brownian motion with local time at 0 denoted by $L$, then, for any positive integrable function $f$ and any bounded adapted (with respect to the filtration of $|B|$) process  $(F_t)$, 
\begin{align}
\lim_{t'\rightarrow +\infty} \frac{\bP[F_t f(L_{t'})]}{\bP[f(L_{t'})]} 
= \bP \sbra{ F_t \frac{M_t}{M_0} } =: \bQ[F_t] 
\label{eq: penal}
\end{align}
where $M$ is the classical Az\'ema-Yor martingale (see \cite{AY}):
\begin{align}
M_t = f(L_t) |B_t| + \int_{0}^{+\infty} f(x+L_t) \d x,\qquad t\geq0.
\nonumber
\end{align}
Then, under the new penalized measure $\bQ$, the random variable $L_\infty$ is finite and the penalized process is seen to be transient. In fact, its paths may roughly be described as the concatenation of a weighted reflected Brownian bridge and a three-dimensional Bessel process; see \cite{MR2253307}.

This classical example on the reflected Brownian motion 
was generalized to many other processes: 
we may refer in particular to Debs \cite{MR2599215} for random walks, 
Najnudel, Roynette and Yor \cite{NRY} for Markov chains and Bessel processes, 
Yano, Yano and Yor \cite{MR2552915} for stable processes, 
or Salminen and Vallois \cite{MR2540855} and Profeta \cite{MR2968676, CP2} 
for linear diffusions. 
In most of these papers, a (sometimes implicit but rather strong) condition is made on the considered process, basically stating that a given quantity is regularly varying. We shall give a short account of this classical setup at the end of this introduction.

In this paper, we shall focus on local time penalizations for generalized linear diffusions, but without any assumption of regular variations. To do so, we shall replace the constant time $t'$ 
by a {\em clock} 
$\tau = \tau_{\lambda} $, 
i.e., a family of random times 
parametrized by $ \lambda $ in a directed set such that 
$ \tau = \tau_{\lambda} $ tends to $+\infty$ a.s. 
We will deal, for example, with the {\em exponential clock} $ \tau = \ve_q := \ve/q $, $ q>0 $ 
for an independent standard exponential random variable $ \ve $ 
and with a parameter $ q>0 $ equipped with the decreasing direction; 
in fact, $ \ve_q \to + \infty $ a.s. as $ q \down 0 $. 
We shall therefore 
study the limit such as 
\begin{align}
\lim_{\tau \to + \infty } \frac{\bP[F_t f(L_{\tau})]}{\bP[f(L_{\tau})]} 
\nonumber
\end{align}
(the limit here is taken along $ \tau = \tau_{\lambda} $ 
with respect to the parameter $ \lambda $, 
and it can also be understood in the sense of Moore--Smith convergence) 
and prove in particular that 
\begin{itemize}
\item[] \vspace{-1cm}
\item[$i)$] no conditions are needed on the characteristics of the diffusion for the existence of a limit, 
\item[$ii)$] the limit depends, in general, of the chosen clock.
\end{itemize}
Examples of such results already appear in the literature, 
dealing with processes conditioned to avoid 0, 
i.e. with the function $f(u)=1_{\{u=0\}}$. 
We refer to Knight \cite{MR0253424} for Brownian motions, 
Chaumont and Doney \cite{MR2164035} and Doney \cite{MR2164035} for L\'evy processes, 
and Yano and Yano \cite{YYrenorm} for diffusions. 

\subsection{Notations}
We start with some notations which are borrowed from \cite{YYrenorm}. Let us consider a generalized one-dimensional diffusion $X$ 
(in the sense of Watanabe \cite{MR0365731}) defined on an interval $I$ 
whose left boundary is 0, with scale function $s(x)=x$ and speed measure $\d m(x)$. 
We assume that the function $m:[0,+\infty) \to [0,+\infty]$ 
is non-decreasing, right-continuous, null at 0 and such that 
\begin{align}
m \text{ is }
\begin{cases}
\text{ strictly increasing on }[0, \ell^\prime),\\
\text{ flat and finite on }[\ell^\prime, \ell),\\
\text{ infinite on }[\ell, +\infty) 
\end{cases}
\nonumber
\end{align}
where $0< \ell^\prime\leq \ell\leq +\infty$. 
The choice of the right boundary point of $I$ will depend on $m$; 
see Section \ref{app} for a sum-up of the different situations (as well as some examples), or \cite[Section 2]{YYrenorm} for a detailed explanation. As for the left boundary point, we assume that $0$ is regular-reflecting for $X$. 
This implies in particular that $X$ admits a local time at 0, 
which we shall denote $(L_t, t\geq0)$, 
normalized so that 
\begin{align}
\bP_x \sbra{ \int_0^{\infty } \e^{-qt} \d L_t } 
= r_q(x,0) 
, \quad x \in I', 
\label{eq: Revuz}
\end{align}
where $ r_q(x,y) $ denotes the resolvent density of $ X $ with respect to $ \d m(y) $ 
and where $ I' $ is defined in Section \ref{app} (see also \cite[Section 2]{YYrenorm}). 
From this formula, we easily obtain the Laplace transform of the first hitting time of the level 0 by $X$:
\begin{align}
\bP_x[\e^{-qT_0}] = \frac{r_q(x,0)}{r_q(0,0)} 
, \quad x \in I' , 
\label{eq: Revuz2}
\end{align}
where $ T_0 = \inf \{ t \ge 0 : X_t = 0 \} $. 
The right-continuous inverse $ \eta^0_u = \inf \{ t \ge 0 : L_t > u \} $ of $ L $ is, 
when considered under $ \bP_0 $, a subordinator. 
We can compute its Laplace exponent
easily by using \eqref{eq: Revuz}, as 
\begin{align}
\bP_x \sbra{ \e^{-q \eta^0_u} } = \e^{- u / H(q)} 
\quad \text{with $ H(q) :=  r_q(0,0) $}. 
\label{eq: Revuz3}
\end{align}

Let $\phi_q$ and $\psi_q$ be the two classical eigenfunctions associated to $X$ 
via the integral equations, for $x\in [0,\ell)$: 
\begin{align}
\phi_q(x) =& 1+q \int_0^x dy \int_{(0,y]}\phi_q(z) \d m(z) , 
\label{phiq} \\
\psi_q(x) =& x+q \int_0^x dy \int_{(0,y]}\psi_q(z) \d m(z). 
\nonumber
\end{align}
With these notations, the resolvent density of $X$ is given by 
\begin{align}
r_q(x,y) = r_q(y,x) = H(q) \phi_q(x) \rbra{ \phi_q(y) - \frac{\psi_q(y)}{H(q)} }, 
\quad 0\leq x\leq y,\quad x,y\in I^\prime,
\nonumber
\end{align}
We denote $\displaystyle  m(\infty ) := \lim_{x \to +\infty } m(x) $, $ \pi_0 := 1/m(\infty ) $ and 
\begin{align}
h_q(x) = r_q(0,0)-r_q(0,x) . 
\nonumber
\end{align}
We finally define, following \cite{YYrenorm}, 
\begin{align}
h_0(x) = \lim_{q\downarrow 0} h_q(x) = x- \pi_0\int_0^x m(y) \d y, 
\label{eq:h0}
\end{align}
and we call $ h_0 $ the {\em normalized zero resolvent}.

\subsection{Main results}

We now outline the main results of the paper. 
For simplicity in this introduction, we assume that $ \ell' = \ell = \infty $  
and we take up the following three cases: 
\begin{center}
\begin{tabular}{|l|l|l|}
\hline
Boundary $ \infty $ & $ m(\infty ) $ & $ \int_{(1,\infty )} x \d m(x) $ 
\\ \hline
(i) type-1-natural & $ =\infty $ & $ =\infty $
\\ \hline
(ii) type-2-natural & $ <\infty $ & $ =\infty $
\\ \hline
(iii) entrance & $ <\infty $ & $ <\infty $
\\ \hline
\end{tabular}
\end{center}
For instance, the Brownian motion reflected at 0 
is an example of case (i), where $ \pi_0=0 $ and $ h_0(x)=x $. 

Let $ \cF^0_t = \sigma(X_s: s \le t) $ and $ \cF_t = \cF^0_{t+} $. 
Let $ \mathcal{L}^1_+ $ denote the set of non-negative functions $ f $ on $ [0,\infty ) $ such that $ \int_0^{\infty } f(u) \d u < \infty $. 
For each choice of a clock 
$ \tau = \tau_{\lambda} $ parametrized by $ \lambda $ 
and a function $ f \in \mathcal{L}^1_+ $, 
we shall find 
a positive function $ \rho(\lambda) $, 
a supermartingale $ N^{f} $ and a martingale $ M^{f} $ with respect to $(\cF_t)$ 
such that the following convergences hold: 
for any $ t>0 $ and any bounded adapted process $ (F_t) $, we have 
\begin{align*}
\rho(\lambda) \bP_x[F_t f(L_{\tau_{\lambda}});t<\tau_{\lambda}] 
\tend{}{\tau = \tau_{\lambda} \to \infty } 
\bP_x[F_t N^{f}_t] 
\quad \text{and} \quad 
\rho(\lambda) \bP_x[F_t f(L_{\tau_{\lambda}})] 
\tend{}{\tau = \tau_{\lambda} \to \infty } 
\bP_x[F_t M^{f}_t] . 
\end{align*}
From these formulae we can obtain the penalization limits of the form 
\begin{align}
\frac{\bP_x[F_t f(L_{\tau});t<{\tau}]}{\bP_x[f(L_{\tau})]} 
\tend{}{\tau \to \infty } 
\bP_x \sbra{ F_t \frac{N^{f}_t}{N^{f}_0} } 
\quad \text{and} \quad 
\frac{\bP_x[F_t f(L_{\tau})]}{\bP_x[f(L_{\tau})]} 
\tend{}{\tau \to \infty } 
\bP_x \sbra{ F_t \frac{M^{f}_t}{M^{f}_0} }  
\nonumber
\end{align}
which allow to construct new (sub-)probabilities.
Note that $ N^{f}_t $ and $ M^{f}_t $ may differ according to 
a particular choice of a clock. 

We shall study such penalization limits with four different clocks as below.

\noindent
1$ ^{\circ} $) \textbf{The exponential clock (see Theorem \ref{2}).}\\
Recall that $\ve$ denotes a standard exponential random variable 
which is independent of the considered diffusion $ X $ 
and $ \ve_q = \ve/q $ with a parameter $q>0$. 
We may adopt $ \{ \ve_q : q>0 \} $ as a clock 
since $ \ve_q \to \infty $ a.s. as $ q \down 0 $. 

Here we assume for notational simplicity 
that $ \ve $ is defined on the same probability space as $ X $; 
in particular, the expectations are taken on $X$ and $ \ve_q = \ve/q $ at the same time. 
Note that $ \ve_q $ is independent of $ \sigma(\bigcup_t \cF_t) $ 
because the filtration $ (\cF_t) $ is generated by $ X $.

\begin{Thm}\label{theo:I1}
Let $f\in \mathcal{L}^1_+$ and $x \ge 0$. 
For any $t>0$ and any bounded adapted process $(F_t)$, 
\begin{align}
H(q) \bP_x \sbra{ F_t f(L_{\ve_q}) ; t<\ve_q } 
\tend{}{q\downarrow 0} 
\bP_x\left[F_t N_t^{h_0, f}\right] 
\quad \text{and} \quad 
H(q) \bP_x \sbra{ F_t f(L_{\ve_q}) } 
\tend{}{q\downarrow 0} 
\bP_x\left[F_t M_t^{h_0, f}\right] 
\nonumber
\end{align}
where the $\bP_x$-supermartingale $N^{h_0, f}$ 
and the $\bP_x$-martingale $M^{h_0, f}$ are defined by 
\begin{align}
N_t^{h_0, f} = h_0(X_t)f(L_t)+ \int_0^{+\infty} f(L_t+u) \d u, 
\qquad t\geq0 
\label{eq:N}
\end{align}
and
\begin{align}
M_t^{h_0,f}= N_t^{h_0, f} + \pi_0 \int_0^t f(L_u) \d u, \qquad t\geq0.
\label{eq:M}
\end{align}
\end{Thm}

\noindent
2$ ^{\circ} $) \textbf{The hitting time clock (see Theorem \ref{2-1}).}\\
For $a\in I$, let $T_a=\inf\{t\geq0: X_t=a\}$ denote 
the first hitting time of $a$ by $X$. 
We may adopt $ \{ T_a: a \ge 0 \} $ as a clock 
since $ T_a \to \infty $ a.s. as $ a \to \infty $. 

\begin{Thm} \label{theo:I2}
Assume that $\infty$ is natural. 
Let $f\in \mathcal{L}^1_+$ and $x \ge 0$. 
For any $t>0$ and any bounded adapted process $(F_t)$, 
\begin{align}
a \bP_x \sbra{ F_t f(L_{T_a});t<T_a } 
\tend{}{a\uparrow +\infty} 
\bP_x\left[F_t M_t^{s, f}\right] 
\quad \text{and} \quad 
a \bP_x \sbra{ F_t f(L_{T_a}) } 
\tend{}{a\uparrow +\infty} 
\bP_x\left[F_t M_t^{s, f}\right] 
\nonumber
\end{align}
where $M^{s, f}$ is the $\bP_x$-martingale defined by 
\begin{align}
M_t^{s, f} = X_t f(L_t) + \int_0^{+\infty} f(L_t+u) \d u, 
\qquad t\geq0.
\label{eq:Ms}
\end{align}
(Here by the superscript $ s $ we mean the scale function $ s(x)=x $.) 
\end{Thm}

\noindent
3$ ^{\circ} $) \textbf{The inverse local time clock (see Theorems \ref{theo:etaa} and \ref{ulim2}).}\\
For $a \ge 0$, let $(L_t^a, t\geq0)$ denote the local time of $X$ at level $a$, 
and define its right-continuous inverse: 
\begin{align}
\eta_u^a = \inf\{t\geq0,\; L_t^a >u\}. 
\nonumber
\end{align}
We may adopt as a clock $ \{ \eta_u^a: a \ge 0 \} $ for a fixed $ u>0 $, 
since $ \eta_u^a \ge T_a \to \infty $ a.s. as $ a \to \infty $. 

\begin{Thm} \label{theo:I3}
Assume that $\infty $ is natural. 
Let $ f \in \mathcal{L}^1_+ $, $ x \ge 0 $ and $u>0$. 
For any $t>0$ and any bounded adapted process $(F_t)$, 
\begin{align}
a \bP_x \sbra{ F_t f(L_{\eta^a_u}); t< \eta_u^a } 
\tend{}{a\uparrow +\infty} 
\bP_x\left[F_t M_t^{s, f}\right] 
\quad \text{and} \quad 
a \bP_x \sbra{ F_t f(L_{\eta^a_u}) } 
\tend{}{a\uparrow +\infty} 
\bP_x\left[F_t M_t^{s, f}\right] 
\nonumber
\end{align}
where $M^{s, f}$ is the $\bP_x$-martingale defined above. 
\end{Thm}

We obtain here the same penalization limit as that of Theorem \ref{theo:I2}. 
In spite of this fact, 
the proofs of the two theorems are quite different; 
we know that $\eta_u^{a} \law T_a + \widetilde{\eta}_u^a $ 
where $\widetilde{\eta}_u^a$ is the inverse local time at $ a $ 
of (an independent copy of) $X$ started at $a$ 
but this fact cannot reduce Theorem \ref{theo:I3} to Theorem \ref{theo:I2}.

We may also adopt as a clock $ \{ \eta_u^a: u \ge 0 \} $ for a fixed $ a>0 $, 
since $ \eta_u^a \to \infty $ a.s. as $ u \to \infty $. 
For this clock we only consider the weights $ f(L_{\eta^a_u}) $ 
for $ f(u) = \e^{- \beta u} $. 

\begin{Thm} \label{theo:I4}
Let $x, a\ge 0$ and $\beta>0$. 
For $t>0$ and any bounded adapted process $(F_t)$, 
\begin{align}
\e^{\frac{\beta u}{1 + a \beta}} 
\bP_x \sbra{ F_t \e^{- \beta L_{\eta^a_u}}; t<\eta^a_u } 
\tend{}{u\uparrow +\infty} 
\bP_x[F_t M_t^{\beta, a}] 
\quad \text{and} \quad 
\e^{\frac{\beta u}{1+a \beta}} 
\bP_x\sbra{ F_t \e^{- \beta L_{\eta^a_u}} } 
\tend{}{u\uparrow +\infty} 
\bP_x[F_t M_t^{\beta, a}] 
\nonumber
\end{align}
where $M^{\beta, a}$ is the $\bP_x$-martingale defined by
\begin{align}
M_t^{\beta, a} =\frac{1+\beta(X_t\wedge a)}{1+\beta a} 
\exp \rbra{ -\beta L_t + \frac{\beta}{1+\beta a} L_t^a } , 
\qquad t\geq0. 
\nonumber
\end{align}
\end{Thm}

\subsection{Comparison among the penalized measures}

Let $ f \in \mathcal{L}^1_+ $. 
We denote by $ \bQ^{h_0,f}_x $ the probability measure such that 
for any $ t>0 $, 
\begin{align}
\bQ^{h_0,f}_x(A) = \bP_x \sbra{ 1_A \frac{M^{h_0,f}_t}{M^{h_0,f}_0} } 
, \quad A \in \cF_t . 
\nonumber
\end{align}
We also define 
$ \bQ^{s,f}_x $ (resp. $ \bQ^{\beta,a}_x $) 
by replacing $ M^{h_0,f} $ by $ M^{s,f} $ (resp. $ M^{\beta,a} $); 
when we speak of $ M^{s,f} $ or $ \bQ^{s,f}_x $ we always assume $ \infty $ is natural. 
Let us compare these three measures. 
We assume for simplicity that $ \int_0^{\infty } f(u) \d u = 1 $ and 
we notice that $ M^{h,f}_0 = N^{h,f}_0 = 1 $ for  both $ h=h_0 $ and $ h=s $. 

When $ \pi_0=0 $, we have $ h_0 = s $ 
and hence the processes $ M^{h_0,f} $ and $ M^{s,f} $ agree, 
which implies $ \bQ^{h_0,f}_x = \bQ^{s,f}_x $. 
Therefore, in this specific situation, Theorem 1.1, 1.2 and 1.3 
yield the same penalized process. 
This is for instance the case for the reflected Brownian motion. 
When $ \pi_0>0 $, i.e., $ 0 $ is positive recurrent, 
the two processes  $ M^{h_0,f} $ and $ M^{s,f} $ disagree; they can be considered 
to be different generalizations of the Az\'ema-Yor martingales. 

Let us focus on the laws of the total local time $L_\infty$. 
The aim of the penalization procedure was to reduce the local time of the original process, so we shall check, for each case, if the strength of the penalization was strong enough to make the total local time of the penalized process finite. For simplicity we only consider the case $ x=0 $. 
For $ u \ge 0 $, from the optional stopping theorem, 
\begin{align}
\bQ_0^{h_0,f}(L_t\geq u) 
= \bP_0 \sbra{ M_t^{h_0,f} 1_{\{L_t \geq u\}} } 
= \bP_0 \sbra{ M_{\eta_u^0}^{h_0,f} 1_{\{\eta_u^0\leq t\}} } , 
\label{bQ_0^h_0,f(L_tu)}
\end{align}
and then from the monotone convergence theorem, 
\begin{align}
\bQ_0^{h_0,f}(L_{\infty } \ge u) 
= \lim_{t \to \infty } \bQ_0^{h_0,f}(L_t \ge u) 
= \bP_0 \sbra{ M_{\eta_u^0}^{h_0,f} } . 
\label{bQ_0^h_0,f(L_iu)}
\end{align}
By definition of $ M^{h_0,f} $, since $ X_{\eta^0_u} = 0 $ and since $ L_{\eta^0_u} = u $, 
\begin{align}
\bP_0 \sbra{ M_{\eta_u^0}^{h_0,f} } 
= \int_u^{+\infty} f(y) \d y 
+ \pi_0 \, \bP_0 \sbra{ \int_0^{\eta_u^0} f(L_r) \d r } . 
\nonumber 
\end{align}
By the change of variables $ r=\eta^0_s $ and since $\bP_0[\eta_u^{0}] = m(\infty) u$, 
\begin{align}
\pi_0 \, \bP_0 \sbra{ \int_0^{\eta_u^0} f(L_r) \d r } 
= \pi_0 \, \bP_0 \sbra{ \int_0^u f(s) \d \eta^0_s } 
= \pi_0 m(\infty ) \int_0^u f(s) \d s . 
\nonumber
\end{align}
We therefore deduce that 
\begin{align}
\begin{cases}
\bQ^{h_0,f}_0(L_{\infty } \in \d u) = f(u) \d u & (\text{if $ \pi_0=0 $}), \\
\bQ^{h_0,f}_0(L_{\infty } = \infty ) = 1 & (\text{if $ \pi_0>0 $}). 
\end{cases}
\nonumber
\end{align}
In other words, when $\pi_0=0$, we went from an original process spending an infinite amount of time at 0 to a penalized process whose total local time at 0 is finite.
A similar analysis shows that:
\begin{align}
\bQ^{s,f}_0(L_{\infty } \in \d u) = f(u) \d u 
\quad \text{(in any case $ \pi_0=0 $ or $ \pi_0>0 $)} . 
\nonumber
\end{align}
By the same argument as \eqref{bQ_0^h_0,f(L_tu)}-\eqref{bQ_0^h_0,f(L_iu)}, we can obtain 
\begin{align}
\bQ^{\beta,a}_0(L_{\infty } \ge u) 
= \bP_0 \sbra{ \frac{M^{\beta,a}_{\eta^0_u}}{M^{\beta,a}_0} } 
= \exp \rbra{ -\beta u } \bP_0 \sbra{ \exp \rbra{ \frac{\beta}{1 + \beta a} L^a_{\eta^0_u} }  } . 
\nonumber
\end{align}
Then, by a direct adaptation of the proof of Lemma \ref{LawLeta} 
and by analytic continuation, we obtain 
\begin{align}
\bP_0 \sbra{ \exp \rbra{\lambda L^a_{\eta^0_u} } } 
= \exp \rbra{\frac{\lambda u}{1-\lambda a}} 
, \quad \lambda < 1/a . 
\nonumber
\end{align}
Hence we deduce that $ \bQ^{\beta,a}_0(L_{\infty } \ge u) = 1 $ for all $ u>0 $, 
from which we conclude that 
\begin{align}
\bQ^{\beta,a}_0(L_{\infty } = \infty ) = 1 
\quad \text{(in any case $ \pi_0=0 $ or $ \pi_0>0 $)} . 
\nonumber
\end{align}
\noindent
Finally, we see that the three cases we have studied provide three different behaviors: 
\begin{enumerate}
\item[]\vspace{-1cm}
\item $\bQ^{s,f}_0$ is a strong penalization: $\bQ^{s,f}_0(L_\infty<\infty)=1$ whatever the case.
\item $\bQ^{h_0,f}_0$ is an intermediate penalization: $\bQ^{h_0,f}_0(L_\infty<\infty)=1 \text{ or } 0$ according to $\pi_0$.
\item $ \bQ^{\beta,a}_0$ is a weak penalization: $\bQ^{\beta,a}_0(L_{\infty } <\infty) = 0$ whatever the case.
\end{enumerate}

\subsection{The case of a constant clock}

One of the most interesting features of the previous theorems is their universality: the results we obtain are the same for any diffusions, and no conditions are imposed on their characteristics. Such universality is no longer true for constant clocks.

\noindent
Assume for instance that the diffusion $X$ is recurrent, with $\infty$ a natural boundary. Then, to get the convergence
\begin{equation}\label{eq:const}
\frac{\bP_x[F_t f(L_{t'})]}{\bP_x[f(L_{t'})]} 
\tend{}{t' \to \infty } 
\bP_x \sbra{ F_t \frac{M^{s,f}_t}{M^{s,f}_0} } 
\end{equation}
like in the reflected Brownian case (with $M^{s,f}$ given by (\ref{eq:Ms})), one needs to add an extra assumption. A sufficient condition is given by Salminen--Vallois \cite{MR2540855}, who require that 
\begin{align}
\text{the (normalized) L\'evy measure $\nu$ of the subordinator $\eta^0$  be subexponential.} 
\label{eq: SVcond}
\end{align}
It may be shown that (\ref{eq: SVcond}) holds for instance if the function $1/(qH(q))$ is regularly varying at $q\rightarrow 0^+$ with exponent in $(0,1)$. We refer here to Profeta \cite{MR2968676}, where other conditions are  also discussed.

\noindent
If such extra assumptions are not fulfilled, then generally the martingale we obtained is different and the asymptotics of $\bP_x[f(L_{t'})]$ depends on the function $f$ (as in Theorem \ref{theo:I4}). We give below two examples.
\begin{enumerate}
\item Take to simplify $f(u)=1_{\{u=0\}}$. Then, under Salminen--Vallois' condition \eqref{eq: SVcond}, 
$$\bP_x(T_0>t')\mathop{\sim}\limits_{t' \to \infty } x\, \nu((t',+\infty))\qquad \text{ and }\qquad
\frac{\bP_x[F_t1_{\{T_0>t'\}}]}{\bP_x(T_0>t')} 
\tend{}{t' \to \infty } 
\bP_x \sbra{F_t\frac{X_t}{x} 1_{\{T_0>t\}}} $$
which is a special instance of (\ref{eq:const}).
If we now consider an Ornstein--Uhlenbeck process $Z$ with parameter $\gamma>0$, for which their condition \eqref{eq: SVcond} on $\nu$ does not hold since for $\varepsilon$ small enough
$$\lim_{x\rightarrow +\infty} e^{\varepsilon x}\nu(x,+\infty) =0  \qquad (\text{instead of } +\infty),$$ then 
$$\bP_x(T_0>t')\mathop{\sim}\limits_{t' \to \infty } 2x\sqrt{\frac{\gamma }{\pi}}e^{-\gamma t'}\qquad \text{ and }\qquad\frac{\bP_x[F_t1_{\{T_0>t'\}}]}{\bP_x(T_0>t')} 
\tend{}{t' \to \infty } 
\bP_x \sbra{F_t\frac{Z_t}{x} e^{\gamma t}\,1_{\{T_0>t\}}},$$
in which case the structure of the martingale is different.
\item Assume now that $X$ is a Brownian motion reflected on $[0,1]$ and take $f(u)=e^{-\alpha u}$ with $\alpha>0$. Then 
$$\bP_x[ e^{-\alpha L_{t'}}] \mathop{\sim}\limits_{t' \to \infty } \kappa \,e^{-r^2 t'} \cos(r\sqrt{2}(1-x))$$
where $r$ is defined as the unique solution in $(0, \pi/(2\sqrt{2}))$ of the equation $\alpha = r\sqrt{2}\tan(r\sqrt{2})$, and $\kappa$ is an explicit constant. In this case, the asymptotics depend on the chosen weight, i.e. on $\alpha$ here.
This yields thus a martingale different from (\ref{eq:Ms}) 
$$\frac{\bP_x[F_t e^{-\alpha L_{t'}}]}{\bP_x[ e^{-\alpha L_{t'}}]} 
\tend{}{t' \to \infty } 
\bP_x \sbra{ F_t  \, e^{r^2 t-\alpha L_t} \frac{\cos\left(r\sqrt{2}(1-X_t)\right)}{\cos\left(r\sqrt{2}(1-x)\right)}}.$$
Such a result may be generalized to other reflected diffusions on $[0,1]$, under the assumption that the analytic continuation of $H(q)$  is smaller (at infinity) than a negative power of $|q|$ on a given strip on the complex plane, see Profeta \cite{CP2}.
 \end{enumerate}

\subsection{Organization}

The remaining of this paper is organized as follows. 
The local time penalizations are studied 
with an independent exponential clock in Section \ref{sec: exp}, 
then with a hitting time clock in Section \ref{sec: ht} 
and finally with inverse local time clocks in Section \ref{sec: ilt}. 
In Section \ref{uni}, we discuss some features of the penalized processes, and give a path decomposition via some universal $ \sigma $-finite measures. 
In Section \ref{ex}, we characterize the limit measure for an exponential weight, in which case the penalized process remains a generalized diffusion. The final section, Section \ref{app}, 
is an appendix on our boundary classification, with some examples.

\section{Local time penalization with an exponential clock} \label{sec: exp} 

The proof of Theorem \ref{theo:I1} will be decomposed in three parts. We shall first study the law of $L_{\ve_q}$, then obtain some a.s. convergence results, and finally extend them to get the $\mathcal{L}^1$ convergence.
Once the law of $ L_{\ve_q} $ has been computed via the excursion theory, 
the method consists in decomposing $ \bP_x[f(L_{\ve_q})|\cF_t] $ according to 
$ \{ t<\ve_q \} $ or $ \{ t \ge \ve_q \} $, 
and then studying both limits separately.

 \subsection{The law of $L_{\ve_q}$}
 
 We start by computing $\bP_x[f(L_{\ve_q})]$.

\begin{Lem} \label{lem1}
Let $ f $ be a non-negative measurable function. 
Let $ q>0 $ and $ x \in I $. Then 
\begin{align}
\bP_x[f(L_{\ve_q})] 
= \frac{1}{H(q)} \cbra{ 
h_q(x) f(0) + \frac{r_q(x,0)}{r_q(0,0)} \int_0^{\infty } \e^{-u/H(q)} f(u) \d u } . 
\label{eq: 1}
\end{align}
\end{Lem}

\Proof{
Using the excursion theory, we have, when starting the diffusion at 0: 
\begin{align}
\bP_0 \sbra{ \int_0^{\infty } f(L_t) q \e^{-qt} \d t } 
=& \bP_0 \sbra{ \sum_{u} \int_{\eta^0_{u-}}^{\eta^0_u} f(u) q \e^{-qt} \d t } 
\nonumber \\
=& \bP_0 \sbra{ \sum_{u} f(u) \e^{- q \eta^0_{u-}}\left(1- \e^{-q T_0(p(u))}\right) } 
\nonumber 
\end{align}
where $ \eta^0_u $ denotes the inverse local time defined around \eqref{eq: Revuz3}, 
$ p $ the excursion point process 
and $T_0(p(s)) = \eta^0_u - \eta^0_{u-} $ the length of the excursion  of $p(s)$. 
Now, denoting by $\vn$ the It\^o's excursion measure, we have 
\begin{align}
\vn \sbra{ 1 - \e^{-q T_0} } = \frac{1}{H(q)} , 
\label{eq: Revuz5}
\end{align}
which follows from \eqref{eq: Revuz3} 
and the formula $ \eta^0_u = \sum_{s \le u} T_0(p(s))) $. 
Using the Master Formula (see \cite[p.475]{RY}), we further obtain:
\begin{align}
\bP_0 \sbra{ \int_0^{\infty } f(L_t) q \e^{-qt} \d t } =& \bP_0 \sbra{ \int_0^{\infty } f(u) \e^{- q \eta^0_u} \d u } 
\vn \sbra{ 1 - \e^{-q T_0} } 
\nonumber \\
=& \int_0^{\infty } f(u) \e^{- u/H(q)} \d u \cdot \frac{1}{H(q)}. 
\nonumber
\end{align}
Here we used \eqref{eq: Revuz3} and \eqref{eq: Revuz5}. 
The Markov property then yields the announced result for any starting point $x\in I$:
\begin{align}
\bP_x[f(L_{\ve_q})] 
=& \bP_x \sbra{ \int_0^{\infty } f(L_t) q \e^{-qt} \d t } 
\nonumber \\
=& \bP_x \sbra{ \int_0^{T_0} f(L_t) q \e^{-qt} \d t } 
+ \bP_x \sbra{ \e^{-q T_0} } \bP_0 \sbra{ \int_0^{\infty } f(L_t) q \e^{-qt} \d t } 
\nonumber \\
=& f(0) \cbra{ 1 - \frac{r_q(x,0)}{r_q(0,0)} } 
+ \frac{r_q(x,0)}{r_q(0,0)} \cdot \int_0^{\infty } f(u) \e^{- u/H(q)} \d u \cdot \frac{1}{H(q)} . 
\nonumber
\end{align}
where we used \eqref{eq: Revuz2} and the fact that $ L_t = 0 $ for $ t \le T_0 $. 
}

Letting $q\downarrow0$, we then deduce the following formulae in the transient and positive recurrent cases. 
(Note that in the null recurrent case, the limit equals $+\infty$).

\begin{Thm} \label{1}
Let $ f $ be a non-negative measurable function and let $ x \in I $. 
Then the following assertions hold: 
\begin{enumerate}
\item 
If $ \ell<\infty $, i.e., 0 is transient, then 
\begin{align}
\bP_x[f(L_{\infty })] = \frac{1}{\ell} \cbra{ x f(0) 
+ \rbra{ 1 - \frac{x}{\ell} } \int_0^{\infty } \e^{-u/\ell} f(u) \d u } . 
\label{eq: 2}
\end{align}
\item 
If $ \pi_0>0 $, i.e., 0 is positive recurrent, then 
\begin{align}
\bP_x \sbra{ \int_0^{\infty } f(L_t) \d t } = \frac{1}{\pi_0} \cbra{ h_0(x) f(0) 
+ \int_0^{\infty } f(u) \d u } . 
\label{eq: 3}
\end{align}
\end{enumerate}
\end{Thm}

\Proof{
(i) We first suppose that $ f $ is bounded and set $ g = \sup \{ t: X_t=0 \} $. 
We see that, almost surely, $ f(L_{\ve_q}) = f(L_g) = f(L_{\infty }) $ 
for $ q>0 $ small enough. 
Thus 
\begin{align}
\bP_x[f(L_{\ve_q})] \tend{}{q \down 0} \bP_x[f(L_{\infty })] 
\nonumber
\end{align}
by the dominated convergence theorem 
(here we do not need continuity of $ f $). 
Equation \eqref{eq: 2} then follows by letting $q\downarrow 0$ in \eqref{eq: 1}. To remove the boundedness assumption, it remains to apply Formula \eqref{eq: 2} with $ f \wedge n $ and then let $ n \to \infty $. 

(ii) 
We may rewrite \eqref{eq: 1} as 
\begin{align}
\bP_x \sbra{ \int_0^{\infty } f(L_t) \e^{-qt} \d t } 
= \frac{1}{qH(q)} \cbra{ 
h_q(x) f(0) + \frac{r_q(x,0)}{r_q(0,0)} \int_0^{\infty } \e^{-u/H(q)} f(u) \d u } . 
\nonumber
\end{align}
Equation \eqref{eq: 3} then follows by letting $ q \down 0 $ and  applying the monotone convergence theorem.
}

\subsection{A.s. convergence for the exponential clock}
Recall that $ \ve_q $ is independent of $ \sigma(\bigcup_t \cF_t) $. 

\begin{Lem} \label{lem2}
Let $ f \in \mathcal{L}^1_+ $ and $ x \in I $. For $ q>0 $, set 
\begin{align}
N^q_t = H(q) \bP_x \sbra{ f(L_{\ve_q}) 1_{\{ t<\ve_q \}} | \cF_t } , 
\qquad 
M^q_t = H(q) \bP_x \sbra{ f(L_{\ve_q}) | \cF_t } 
\nonumber
\end{align}
and set 
\begin{align}
N^{h_0,f}_t 
=& h_0(X_t) f(L_t) + \rbra{ 1-\frac{X_t}{\ell} } \int_0^{\infty } \e^{-u/\ell} f(L_t+u) \d u , 
\label{eq: Nh0ft} \\
M^{h_0,f}_t =& N^{h_0,f}_t + A^{h_0,f}_t , 
\label{eq: 4} \\
A^{h_0,f}_t =& \pi_0 \int_0^t f(L_u) \d u 
\nonumber
\end{align}
(notice that \eqref{eq:N} and \eqref{eq:M} are the special cases 
for $ \ell=\ell'=\infty $ 
of \eqref{eq: Nh0ft} and \eqref{eq: 4}, respectively). 
Then the following assertions hold: 
\begin{enumerate}
\item 
$ N^q_t \to N^{h_0,f}_t $ and $ M^q_t \to M^{h_0,f}_t $, $ \bP_x $-a.s. as $ q \down 0 $; 
\item 
$ (N^{h_0,f}_t) $ is a $ \bP_x $-supermartingale. 
\end{enumerate}
\end{Lem}

\Proof{
In what follows in this section we sometimes write to simplify $ N_t $, $ M_t $ and $ A_t $ 
 for $ N^{h_0,f}_t $, $ M^{h_0,f}_t $ and $ A^{h_0,f}_t $, respectively. 

(i) Since $(L_t)$ is an additive functional and $\ve_q$ an exponential random variable, we have  
\begin{align}
N^q_t 
=& H(q) \e^{-qt} \left. \bP_{X_t}[ f(a + L_{\ve_q}) ] \right|_{a=L_t} 
\nonumber \\
=& \e^{-qt} \cbra{ 
h_q(X_t) f(L_t) + \frac{r_q(X_t,0)}{r_q(0,0)} 
\int_0^{\infty } \e^{-u/H(q)} f(L_t + u) \d u } 
\nonumber
\end{align}
where we have used Lemma \ref{lem1} with  $f(a+\cdot) \in \mathcal{L}^1_+ $. 
It is now clear that $ N^q_t \tend{}{q \down 0} N_t $, $ \bP_x $-a.s. 
Since 
\begin{align}
A^q_t 
:=& M^q_t - N^q_t 
\nonumber \\
=& H(q) \bP_x[f(L_{\ve_q}) 1_{\{ \ve_q \le t \}} | \cF_t ] 
\nonumber \\
=& qH(q) \int_0^t f(L_u) \e^{-qu} \d u , 
\nonumber
\end{align}
we obtain $ A^q_t \tend{}{q \down 0} A_t $ and $ M^q_t \tend{}{q \down 0} M_t $, $ \bP_x $-a.s.

(ii) 
Since for $ s \le t $ we have $ 1_{\{ t<\ve_q \}} \le 1_{\{ s<\ve_q \}} $, 
we easily see that $ (N^q_t) $ is a $ \bP_x $-supermartingale. 
For $ s \le t $, we apply Fatou's lemma to obtain 
\begin{align}
\bP_x[N_t|\cF_s] 
\le \liminf_{q \down 0} \bP_x[N^q_t|\cF_s] 
\le \liminf_{q \down 0} N^q_s = N_s , 
\nonumber
\end{align}
which shows that $ (N_t) $ is a $ \bP_x $-supermartingale. 
}

\subsection{$\mathcal{L}^1$ convergence for the exponential clock}
\begin{Thm} \label{2}
Let $ f \in \mathcal{L}^1_+ $ and $ x \in I $. 
\begin{enumerate}
\item For any finite stopping time $ T $,  there is the $\mathcal{L}^1$ convergence
\begin{align}
N^q_T \tend{}{q \down 0} N^{h_0,f}_T 
\quad \text{in $ \mathcal{L}^1(\bP_x) $} . 
\nonumber
\end{align}
Consequently, 
for any bounded adapted process $ (F_t) $, it holds that 
\begin{align}
\lim_{q \down 0} H(q) \bP_x [F_T f(L_{\ve_q}) ; T<\ve_q] 
= \bP_x[F_T N^{h_0,f}_T] . 
\nonumber
\end{align}
\item Assume furthermore that 
\begin{align}
\bP_x \sbra{ \int_0^T f(L_u) \d u } < \infty . 
\label{eq: Tfin}
\end{align}
Then, we have 
\begin{align}
M^q_T  \tend{}{q \down 0}  M^{h_0,f}_T 
\quad \text{in $ \mathcal{L}^1(\bP_x) $},
\label{MqMh}
\end{align}
and for any bounded adapted process $ (F_t) $, it holds that 
\begin{align}
\lim_{q \down 0} H(q) \bP_x [F_T f(L_{\ve_q})] 
= \bP_x[F_T M^{h_0,f}_T] . 
\nonumber
\end{align}
\item 
Any bounded stopping time satisfies \eqref{eq: Tfin}. In particular, 
\begin{align}
M^{h_0,f}_t 
= h_0(X_t) f(L_t) + \rbra{ 1-\frac{X_t}{\ell} } \int_0^{\infty } \e^{-u/\ell} f(L_t+u) \d u 
+ \pi_0 \int_0^t f(L_u) \d u 
\nonumber
\end{align}
is a $ \bP_x $-martingale and the identity $ N^{h_0,f}= M^{h_0,f} - A^{h_0,f} $ may be regarded as the Doob--Meyer decomposition of the supermartingale $N^{h_0,f} $. 
\end{enumerate}
\end{Thm}

\Proof{(i) Observe first by Fatou's lemma that 
\begin{align}
\bP_x[N_T] 
\le \liminf_{n \to \infty } \bP_x[N_{T \wedge n}] 
\le \bP_x[N_0] < \infty . 
\nonumber
\end{align}
\noindent
Let us compute $ N^q_T $. We have 
\begin{align}
N^q_T 
=& \e^{-qT} h_q(X_T) f(L_T) + \e^{-qT} \frac{r_q(X_T,0)}{r_q(0,0)} 
\int_0^{\infty } \e^{-u/H(q)} f(L_T + u) \d u 
\nonumber \\
=& ({\rm I})_q + ({\rm II})_q . 
\nonumber
\end{align}
We write similarly 
\begin{align}
N_T 
=& h_0(X_T) f(L_T) + \rbra{ 1-\frac{X_T}{\ell} } \int_0^{\infty } \e^{-u/\ell} f(L_T+u) \d u 
\nonumber \\
=& ({\rm I}) + ({\rm II}) . 
\nonumber
\end{align}
\noindent
Since $ ({\rm II})_q \le \int_0^{\infty } f(u) \d u $, 
we may apply the dominated convergence theorem to obtain 
$ ({\rm II})_q \to ({\rm II}) $ in $ \mathcal{L}^1(\bP_x) $. 

\noindent
If $ \pi_0=0 $, then we have 
$ ({\rm I})_q \le X_T f(L_T) = h_0(X_T) f(L_T) \le N_T $. 
If $ \pi_0>0 $ and $ \ell' $ is regular-reflecting, 
then we have $ h_0(x) \ge cx $ with $ c=h_0(\ell')/\ell'>0 $, since $ h_0(x) $ is concave. 
We now have $ ({\rm I})_q \le X_T f(L_T) \le c^{-1} h_0(X_T) f(L_T) \le c^{-1} N_T $. 
In both cases, since $ \bP_x[N_T] \le \bP_x[N_0] < \infty $, 
we may apply the dominated convergence theorem to obtain 
$ ({\rm I})_q \to ({\rm I}) $ in $ \mathcal{L}^1(\bP_x) $. 

\noindent
If $ \pi_0>0 $ and $ \ell' $ is either entrance or natural, 
we have $ ({\rm I})_q \le X_T f(L_T) $. 
Since we see by (ii) of Lemma \ref{lem2} that 
\begin{align}
\bP_x[X_T f(L_T)] 
\le \bP_x[N^{s,f}_T] 
\le \bP_x[N^{s,f}_0] 
= x f(0) + \rbra{ 1 - \frac{x}{\ell} } \int_0^{\infty } \e^{-u/\ell} f(u) \d u < \infty , 
\nonumber
\end{align}
we may apply the dominated convergence theorem to obtain 
$ ({\rm I})_q \to ({\rm I}) $ in $ \mathcal{L}^1(\bP_x) $. 
Therefore we have obtained the former assertion. 
For the latter assertion, we have 
\begin{align}
H(q) \bP_x [F_T f(L_{\ve_q}) ; T<\ve_q] 
= \bP_x [F_T N^q_T] 
\tend{}{q \down 0} 
\bP_x[F_T N_T] . 
\nonumber
\end{align}

(ii) To prove \eqref{MqMh},  it suffices to observe that 
by \eqref{eq: Tfin}, we have 
$ \int_0^T f(L_u) \e^{-qu} \d u \to \int_0^T f(L_u) \d u $ in $ \mathcal{L}^1(\bP_x) $. 
This shows that $ A^q_T \to A_T $ in $ \mathcal{L}^1(\bP_x) $, 
which implies $ M^q_T \to M_T $ in $ \mathcal{L}^1(\bP_x) $. 

(iii) 
Since 
\begin{align}
q^2 \int_0^{\infty } \bP_x \sbra{ \int_0^t f(L_u) \d u } \e^{-qt} \d t 
= \bP_x \sbra{ \int_0^{\infty } f(L_u) q \e^{-qu} \d u } 
= \bP_x \sbra{ f(L_{\ve_q}) } < \infty 
\nonumber
\end{align}
and since $ t \mapsto \bP_x \sbra{ \int_0^t f(L_u) \d u } $ is increasing, 
we see that $ \bP_x \sbra{ \int_0^t f(L_u) \d u } < \infty $ for all $ t \ge 0 $. 
Thus \eqref{MqMh} holds for all constant times, 
which implies that $ M^{h_0,f} $ is a martingale. 
}

%
%

\begin{Rem}
As mentioned in the Introduction, if $ \ell' $ is type-1-natural, then the identity \eqref{eq: 4} becomes 
\begin{align}
M^{h_0,f}_t = X_t f(L_t) + \int_0^{\infty } f(L_t+u) \d u , 
\nonumber
\end{align}
which is nothing else but the Az\'ema--Yor martingale (\cite{AY}). 
In this sense we may regard the identity \eqref{eq: 4} 
as a generalization of the Az\'ema--Yor martingale. 
Another generalization will be given in Theorem \ref{2-1}. 
\end{Rem}

\begin{Rem}
If we take $ f(u) = 1_{\{ u=0 \}} $, we have 
\begin{align}
M^{h_0,f}_t = h_0(X_t) 1_{\{ T_0>t \}} + \pi_0 (T_0 \wedge t) . 
\nonumber
\end{align}
In particular, from the identity $ \bP_x\left[M^{h_0,f}_0\right] = \bP_x\left[M^{h_0,f}_t\right] $, we obtain 
\begin{align}
h_0(x) = \bP_x[h_0(X_t);T_0>t] + \pi_0 \bP_x[T_0 \wedge t] , 
\nonumber
\end{align}
which is nothing else but the first assertion of Theorem 6.4 of \cite{YYrenorm}. 
\end{Rem}

\section{Local time penalization with a hitting time clock} \label{sec: ht}

In this section we assume that $ \ell$ $( =\ell' $) is either entrance or natural. 
Since any point in $ [0,\ell) $ is accessible but $ \ell $ is not, we have 
\begin{align}
\bP_x \rbra{ T_a \to \infty \ \text{as} \ a \up \ell } = 1 . 
\nonumber
\end{align}

We start by computing the law of $ L_{T_a} $ using the formula \eqref{eq: 2} 
for the stopped process at $ a $. 
We then prove a.s. and $ \cL^1 $ convergence results upon separating the cases 
$ \{ t<T_a \} $ and $ \{ t \ge T_a \} $ and using some estimates 
on the cumulative distribution function of $ T_a $. 

\subsection{A.s. convergence for hitting times clocks}

We start by computing the quantity $\bP_x[f(L_{T_a})]$.

\begin{Lem} \label{lem: 2-1}
Let $ f \in \mathcal{L}^1_+ $ and $ x \in I $. 
Then, for any $ a \in I $ with $ x<a $, 
\begin{align}
\bP_x[f(L_{T_a})] = \frac{1}{a} \cbra{ x f(0) 
+ \rbra{ 1-\frac{x}{a} } \int_0^{\infty } \e^{-u/a} f(u) \d u } . 
\label{eq: 21}
\end{align}
\end{Lem}

\Proof{
Let $ \bP^a_x $ 
denote the law of $ X_{\cdot \wedge T_a} $ under $ \bP_x $. 
Then we have 
\begin{align}
\bP_x[f(L_{T_a})] = \bP^a_x[f(L_{\infty })] . 
\nonumber
\end{align}
Since $ \{ X,\bP^a_x \} $ is a diffusion process on $ [0,a] $ 
where $ a $ is a regular-absorbing boundary, 
we may use (i) of Theorem \ref{1} and obtain \eqref{eq: 21}. 
}

\begin{Lem} \label{lem: 2-2}
Let $ f \in \mathcal{L}^1_+ $ and $ x \in I $. 
For any $ a \in I $ with $ x<a $, set 
\begin{align}
N^a_t = a \bP_x \sbra{ f(L_{T_a}) 1_{\{ t<T_a \}} | \cF_t } , 
\quad 
M^a_t = a \bP_x \sbra{ f(L_{T_a}) | \cF_t } 
\nonumber 
\end{align}
and 
\begin{align}
M^{s,f}_t =& X_t f(L_t) 
+ \rbra{ 1 - \frac{X_t}{\ell} } \int_0^{\infty } \e^{-u/\ell} f(L_t+u) \d u . 
\label{eq: 2-7}
\end{align}
Then the following assertions hold: 
\begin{enumerate}
\item 
$ N^a_t \to M^{s,f}_t $ and $ M^a_t \to M^{s,f}_t $, $ \bP_x $-a.s. as $ a \up \ell $; 
\item 
$ (M^{s,f}_t) $ is a $ \bP_x $-supermartingale and is a local $ \bP_x $-martingale. 
\end{enumerate}
\end{Lem}

\Proof{
In what follows in this section we sometimes write  to simplify $ M_t $ 
 for $ M^{s,f}_t $. 

(i) Since $ f(b+\cdot) \in \mathcal{L}^1_+ $, we have, by Lemma \ref{lem: 2-1}, 
\begin{align}
N^a_t 
=& a \left. \bP_{X_t}[ f(b + L_{T_a}) ] \right|_{b=L_t} 1_{\{ t<T_a \}} 
\nonumber \\
=& \cbra{ X_t f(L_t) 
+ \rbra{ 1-\frac{X_t}{a} } \int_0^{\infty } \e^{-u/a} f(L_t + u) \d u } 1_{\{ t<T_a \}} . 
\nonumber
\end{align}
Using that $ T_a \to \infty $ as $ a \up \ell $, 
we then deduce that $ N^a_t \to M_t $, $ \bP_x $-a.s. 
Set 
\begin{align}
A^a_t = M^a_t - N^a_t = a f(L_{T_a}) 1_{\{ T_a \le t \}} . 
\nonumber
\end{align}
Since $ A^a_t \to 0 $, $ \bP_x $-a.s., 
we further obtain that $ M^a_t \to M_t $, $ \bP_x $-a.s.

(ii) 
In the same way as (ii) of Lemma \ref{lem2}, 
we can see that $ (M_t) $ is a $ \bP_x $-supermartingale. 

\noindent
It is obvious that $ (M^a_t) $ is a $ \bP_x $-martingale. 
Let $ \{ a_n \} $ be a sequence of $ I $ such that $ a_n \up \ell $. 
If we take $ \sigma_n = \inf \{ t: X_t > a_n \} $, 
we have $ A^a_{\sigma_n \wedge t} = a f(L_{T_a}) 1_{\{ T_a \le \sigma_n \wedge t \}} = 0 $ 
for any $ a>a_n $, 
so that we have $ M^a_{\sigma_n \wedge t} \to M_{\sigma_n \wedge t} $ 
in $ \mathcal{L}^1(\bP_x) $ as $ a \up \ell $. 
This shows that $ (M_t) $ is a local $ \bP_x $-martingale. 
}

\subsection{$\mathcal{L}^1$ convergence for hitting times clocks}

\begin{Thm} \label{2-1}
Let $ f \in \mathcal{L}^1_+ $ and $ x \in I $. 
\begin{enumerate}
\item For any finite stopping time $ T $, there is the $\mathcal{L}^1$ convergence 
\begin{align}
N^a_T \tend{}{a \up \ell} M^{s,f}_T 
\quad \text{in $ \mathcal{L}^1(\bP_x) $} . 
\label{eq: 2-2}
\end{align}
Consequently, 
for any bounded adapted process $ (F_t) $, it holds that 
\begin{align}
a \bP_x [F_T f(L_{T_a}) ; T<T_a] 
\tend{}{a \up \ell} \bP_x\left[F_T M^{s,f}_T\right] . 
\nonumber
\end{align}
\item Suppose furthermore that $ \ell $ is natural.  Then, for any $ f \in \mathcal{L}^1_+ $ and for any bounded stopping time $ T $, we have  
\begin{align}
M^a_T \tend{}{a \up \ell} M^{s,f}_T 
\quad \text{in $ \mathcal{L}^1(\bP_x) $} 
\label{eq: 2-4}
\end{align}
and  for any bounded adapted process $ (F_t) $, it holds that 
\begin{align}
a \bP_x \sbra{ F_T f(L_{T_a}) } 
\tend{}{a \up \ell} 
\bP_x \sbra{ F_T M^{s,f}_T } . 
\label{eq: 2-5}
\end{align}
In particular, $M^{s,f}$ is a $ \bP_x $-martingale. 
\end{enumerate}
\end{Thm}

\Proof{We start with Point $(i)$. We have 
\begin{align}
N^a_T 
=& X_T f(L_T) 1_{\{ T<T_a \}} 
+ \rbra{ 1-\frac{X_T}{a} } \int_0^{\infty } \e^{-u/a} f(L_T + u) \d u \,1_{\{ T<T_a \}} , 
\nonumber \\
M_T 
=& X_T f(L_T) + \rbra{ 1-\frac{X_T}{\ell} } \int_0^{\infty } \e^{-u/\ell} f(L_T+u) \d u . 
\nonumber
\end{align}
Since $ N^a_T \le M_T $ and since 
\begin{align}
\bP_x[M_T] 
\le \liminf_{n \to \infty } \bP_x[M_{T \wedge n}] 
\le \bP_x[M_0] < \infty , 
\nonumber
\end{align}
we may apply the dominated convergence theorem to obtain \eqref{eq: 2-2}. 
The remaining assertion is obvious. 

To prove Point $(ii)$, we shall rely on the following Lemma:

\begin{Lem}
Suppose that $ \ell $ is natural. Then 
\begin{align}
a \bP_x(T_a \le t) \tend{}{a \up \ell} 0 
\quad \text{for all $ t \ge 0 $}. 
\label{eq: 2-3}
\end{align}
\end{Lem}

\Proof{
If $ \ell<\infty $, i.e., $ \ell $ is type-3-natural, then \eqref{eq: 2-3} is obvious. 

Suppose $ \ell=\infty $. 
Since the Laplace transform of $T_a$ may be written as a ratio of eigenfunctions, 
we have  
\begin{align}
a \bP_x(T_a \le t) 
\le a \e^t \bP_x\left[\e^{-T_a}\right] 
= \e^t \phi_1(x) \cdot \frac{a}{\phi_1(a)} . 
\nonumber
\end{align}
Going back to the integral equation \eqref{phiq} of $\phi_1$ and using that $\ell=\infty$ is natural, we deduce that 
\begin{align}
\phi_1(a) = 1 + \int_0^a \d x \int_{(0,x]} \phi_1(y) \d m(y) 
\ge \int_0^a \d x \int_{(0,x]} \d m(y) 
\tend{}{a \up \ell} \infty 
\nonumber
\end{align}
and, for $ a>1 $, 
\begin{align}
\phi'_1(a) \ge \int_{(0,a]} \d m(x) \int_0^x \phi_1'(y) \d y 
\ge \phi_1'(1) \int_{(1,a]} \d m(x) \int_1^x \d y 
\tend{}{a \up \ell} \infty . 
\nonumber
\end{align}
Thus, by the l'H\^opital's rule, we obtain 
$ a/\phi_1(a) \to 0 $ as $ a \up \ell = \infty $ which yields \eqref{eq: 2-3}. 
}

%
We now come back to the proof of Point $(ii)$ of Theorem \ref{2-1}.
Suppose that $ f \in \mathcal{L}^1_+ $ is bounded. 
Since $ A^a_T \to 0 $, $ \bP_x $-a.s. 
and since 
\begin{align}
\bP_x[A^a_T] \le a \| f \|_{\infty } \bP_x(T_a \le T) 
\tend{}{a \up \ell} 0 , 
\nonumber
\end{align}
we see that $ A^a_T \to 0 $ in $ \mathcal{L}^1(\bP_x) $. 
Hence we obtain \eqref{eq: 2-4} and \eqref{eq: 2-5} in this special case. 

\noindent
We now see that $ \bP_x\left[M^{s,f}_t\right] = \bP_x\left[M^{s,f}_0\right] $, i.e., 
\begin{align}
\begin{split}
& \bP_x \sbra{ X_t f(L_t) + \rbra{ 1 - \frac{X_t}{\ell} } \int_0^{\infty } \e^{-u/\ell} f(L_t+u) \d u } 
\\
=&\; x f(0) + \rbra{ 1 - \frac{x}{\ell} } \int_0^{\infty } \e^{-u/\ell} f(u) \d u 
\end{split}
\label{eq: 2-8}
\end{align}
holds for all $ t \ge 0 $ and all bounded $ f \in \mathcal{L}^1_+ $. 
By considering $ f \wedge n $, taking $ n \to \infty $ 
and applying the monotone convergence theorem, 
we can drop the boundedness assumption 
and obtain \eqref{eq: 2-8} for all $ t \ge 0 $ and all $ f \in \mathcal{L}^1_+ $. 
By (ii) of Lemma \ref{lem: 2-2}, 
we see, for any $ f \in \mathcal{L}^1_+ $, that 
$ (M^{s,f}_t) $ is a $ \bP_x $-supermartingale with constant expectation, 
which turns out to be a $ \bP_x $-martingale. 

\noindent
Let $ f \in \mathcal{L}^1_+ $. 
Since $ (M^{s,f}_t) $ is a $ \bP_x $-martingale, 
we may apply the optional stopping theorem to see that 
\begin{align}
\bP_x[A^a_T] 
=& \bP_x \sbra{ X_{T_a} f(L_{T_a}) ; T_a \le T } 
\nonumber \\
\le& \bP_x \sbra{ M_{T_a} ; T_a \le T } 
\nonumber \\
=& \bP_x \sbra{ M_T ; T_a \le T } 
\tend{}{a \up \ell} 0 . 
\nonumber
\end{align}
Since $ A^a_T \to 0 $, $ \bP_x $-a.s., 
we see that $ A^a_T \to 0 $ in $ \mathcal{L}^1(\bP_x) $. 
Hence we obtain \eqref{eq: 2-4} and \eqref{eq: 2-5} in the general case. 
}

\begin{Rem}
Suppose $ \ell $ is entrance. 
We claim that $ M=M^{s,f} $ is {\em not} a true $ \bP_x $-martingale. 
Indeed, suppose $ M $ were a $ \bP_x $-martingale. 
On the one hand, we would have 
\begin{align}
\bP_x[M_{\eta^0_u}] 
= \lim_{t \to \infty } \bP_x \sbra{ M_{\eta^0_u \wedge t} 1_{\{ \eta^0_u \le t \}} } 
= M_0 = x f(0) + \int_0^{\infty } f(r) \d r . 
\nonumber
\end{align}
On the other hand, since $ X_{\eta^0_u} = 0 $ and  $ \ell = \infty $, we have 
\begin{align}
\bP_x[M_{\eta^0_u}] = \int_u^{\infty } f(r) \d r 
\tend{}{u \to \infty } 0 , 
\nonumber
\end{align}
which would be a contradiction. 
Note that in the special case $ f(u) = 1_{\{ u=0 \}} $ and $ M^{s,f}_t = X_t 1_{\{ T_0>t \}} $ 
this result has already been obtained in Theorem 6.5 of \cite{YYrenorm} in other words: 
$ s(x) = x $ is {\em not} invariant with respect to the stopped process. 
\end{Rem}

\section{Local time penalization with inverse local time clocks} \label{sec: ilt}

We recall that $\eta_u^a$ denotes the  right-continuous inverse of  $(L_t^a, t\geq0)$:
\begin{align}
\eta_u^a = \inf\{t\geq0,\; L_t^a >u\}. 
\nonumber
\end{align}

The proofs given in this Section follow the same pattern as before. 
The main ingredients are excursion theory away from $ a \neq 0 $ 
and some estimates on modified Bessel functions.

\subsection{Limit as $ a $ tends to infinity with $ u $ being fixed}

Suppose $ \ell'$ $(=\ell=\infty) $ is 
either entrance, type-1-natural or type-2-natural. 
We thus have, for any $ x \in I $ and any $ u>0 $, 
\begin{align}
\bP_x(\eta^a_u<\infty ) = 1 
\quad \text{and} \quad 
\eta^a_u \ge T_a \tend{}{a \to \infty } \infty, \text{$ \bP_x $-a.s.;} 
\nonumber
\end{align}
in fact, 
the former is implied by 
\begin{align}
\bP_x[\e^{- q \eta^a_u}] 
\tend{}{q \down 0} 1 
\nonumber
\end{align}
since
$ r_q(a,a) \tend{}{q \down 0} \infty $ and 
\begin{align}
\bP_x[\e^{- q \eta^a_u}] 
= \bP_x[\e^{-q T_a}] \bP_a[\e^{-q \eta^a_u}] 
= \frac{\phi_q(x)}{\phi_q(a)} \exp \rbra{ - \frac{u}{r_q(a,a)} } . 
\label{eq: Eqeta}
\end{align}

For $ \nu \ge 0 $, we denote by $ I_{\nu}(x) $ the modified Bessel function of the first kind, 
which may be represented as a series expansion formula 
(see e.g. \cite{Leb}, eq. (5.7.1) on page 108) by 
\begin{align}
I_{\nu}(x) 
= \sum_{n=0}^{\infty } \frac{(x/2)^{\nu+2n}}{n! \Gamma(\nu+n+1)} 
, \quad x>0 . 
\label{eq: Inu}
\end{align}
We recall the asymptotic formulae (see e.g. \cite{Leb}, Section 5.16): 
\begin{align}
I_{\nu}(x) 
\mathrel{\mathop{\sim}\limits_{x \down 0}} 
\frac{(x/2)^{\nu}}{\Gamma(1+\nu)} 
, \quad 
I_{\nu}(x) 
\mathrel{\mathop{\sim}\limits_{x \to \infty }} 
\frac{\e^x}{\sqrt{2 \pi x}} . 
\label{eq: asympInu}
\end{align}

\begin{Lem} \label{LawLeta}
Let $ a \in (0,\infty ) $. Then the process $ \{ (L_{\eta^a_u})_{u \ge 0},\bP_a \} $ 
is a compound Poisson process with Laplace transform 
\begin{align}
\bP_a \sbra{ \e^{- \beta L_{\eta^a_u}} } 
= \exp \cbra{ - u \int_0^{\infty } (1-\e^{-\beta s}) \frac{1}{a^2} \e^{-s/a} \d s } 
= \e^{- \frac{u \beta}{1 + \beta a}} . 
\label{eq: LTLeta}
\end{align}
For any $ u>0 $ and $ f \in \mathcal{L}^1_+ $, 
\begin{align}
\bP_a[f(L_{\eta^a_u})] 
=& \e^{-u/a} f(0) + \int_0^{\infty } f(y) \rho^a_u(y) \d y , 
\label{eq: LawLeta}
\end{align}
where 
\begin{align}
\rho^a_u(y) 
= \e^{-(u+y)/a} \frac{\sqrt{u/y}}{a} I_1 \rbra{ \frac{2 \sqrt{uy}}{a} } . 
\nonumber
\end{align}
\end{Lem}

\Proof{
Let $ p^a(v) $ denote the point process of excursions away from $ a $ 
and $ \vn^a $ its excursion measure. 
Since $ L $ increases only on the intervals $ (\eta^a_{v-},\eta^a_v) $, we have 
\begin{align}
L_{\eta^a_u} 
= \sum_{v \le u: \, p^a(v) \in \{ T_0<\infty \}} (L_{\eta^a_v}-L_{\eta^a_{v-}}) 
= \sum_{v \le u: \, p^a(v) \in \{ T_0<\infty \}} L_{T_a}(p^a(v)) . 
\label{eq: Letaat}
\end{align}
Since $ \vn^a(T_0<T_a) = 1/a < \infty $, 
the sum of \eqref{eq: Letaat} is a finite sum, 
and so we see that 
$ \{ (L_{\eta^a_u})_{u \ge 0},\bP_a \} $ is a compound Poisson process 
with L\'evy measure 
\begin{align}
\vn^a(L_{T_a} \in \d s;T_0<T_a) . 
\nonumber
\end{align}
By the strong Markov property of $ \vn^a $ 
(see, e.g., \cite[Theorem III.3.28]{Blu}), we have 
\begin{align}
\vn^a(L_{T_a}>s ; T_0<T_a) 
= \vn^a(T_0<\infty ) \bP_0(L_{T_a}>s) 
= \frac{1}{a} \bP_0(L_{T_a}>s) . 
\nonumber
\end{align}
Let $ \lambda^0_a = \inf \{ v: p^0(v) \in \{ T_a < \infty \} \} $. 
Then we have 
\begin{align}
\bP_0(L_{T_a}>s) 
= \bP_0(T_a>\eta^0_s) 
= \bP_0(\lambda^0_a>s) 
= \e^{-s \vn^0(T_a<\infty )} 
= \e^{-s/a} . 
\nonumber
\end{align}
Thus we obtain \eqref{eq: LTLeta}. 

Let $ \{ S_n \} $ be a process 
with i.i.d. increments $ \bP(S_n-S_{n-1}>s) = \e^{-s/a} $ 
such that $ S_0=0 $ 
and let $ N $ be a Poisson variable with mean $ u/a $ 
which is independent of $ \{ S_n \} $. 
Then we have $ L_{\eta^a_u} \law S_N $, 
and hence 
\begin{align}
\bP_a[f(L_{\eta^a_u})] 
=& \bP(N=0) f(0) 
+ \sum_{n=1}^{\infty } \bP(N=n) \bP[f(S_n)] 
\nonumber \\
=& \e^{-u/a} f(0) 
+ \sum_{n=1}^{\infty } \e^{-u/a} \frac{(u/a)^n}{n!} 
\int_0^{\infty } f(y) \frac{(y/a)^{n-1}}{(n-1)!} \e^{-y/a} \frac{\d y}{a} . 
\nonumber
\end{align}
Thus, using \eqref{eq: Inu}, we obtain \eqref{eq: LawLeta}. 
}

\begin{Lem} \label{LawLeta2}
For $ u>0 $, $ x,a \in I $ and $ f \in \mathcal{L}^1_+ $, it holds that 
\begin{align}
\bP_x[f(L_{\eta^a_u})] 
=& \frac{x \wedge a}{a} \bP_a[f(L_{\eta^a_u})] 
+ \rbra{ 1 - \frac{x}{a} }_+ \bP_a[f(\ve_{1/a} + L_{\eta^a_u})] 
\label{eq: LawLeta1-2} \\
=& \frac{x \wedge a}{a} \bP_a[f(L_{\eta^a_u})] 
+ \frac{1}{a} \rbra{ 1 - \frac{x}{a} }_+ \int_0^{\infty } f(y) \tilde{\rho}^a_u(y) \d y , 
\label{eq: LawLeta2}
\end{align}
where 
\begin{align}
\tilde{\rho}^a_u(y) 
= \e^{-(u+y)/a} I_0 \rbra{ \frac{2 \sqrt{uy}}{a} } . 
\nonumber
\end{align}
\end{Lem}

\Proof{
When $ a \le x $, 
we have 
$ \bP_x[f(L_{\eta^a_u})] = \bP_x[f(L_{T_a} + L_{\eta^a_u} \circ \theta_{T_a})] 
= \bP_a[f(L_{\eta^a_u})] $, 
which proves identity \eqref{eq: LawLeta1-2}. 

Suppose $ x < a $. 
Using Lemma \ref{lem: 2-1}, we have 
\begin{align}
\bP_x[f(L_{\eta^a_u})] 
=& \bP_x \sbra{ f \rbra{ L_{T_a} + L_{\eta^a_u} \circ \theta_{T_a} } } 
\nonumber \\
=& \frac{x}{a} \bP_a[f(L_{\eta^a_u})] 
+ \frac{1}{a} \rbra{ 1 - \frac{x}{a} } \bP_a \sbra{ 
\int_0^{\infty } \e^{-v/a} f(v+L_{\eta^a_u}) \d v } , 
\nonumber
\end{align}
which coincides with \eqref{eq: LawLeta1-2}. 
Using the same notation as that of the proof of Lemma \ref{LawLeta}, 
we obtain 
\begin{align}
\bP_a[f(\ve_{1/a} + L_{\eta^a_u})] 
=& \sum_{n=0}^{\infty } \bP(N=n) \bP[f(S_{n+1})] 
\nonumber \\
=& \sum_{n=0}^{\infty } \e^{-u/a} \frac{(u/a)^n}{n!} 
\int_0^{\infty } f(y) \frac{(y/a)^n}{n!} \e^{-y/a} \frac{\d y}{a} . 
\nonumber
\end{align}
Thus, using \eqref{eq: Inu}, we obtain \eqref{eq: LawLeta2}. 
}

By \eqref{eq: asympInu}, there exists a constant $ C $ such that 
\begin{align}
\begin{cases}
I_{\nu}(x) \le C x^{\nu} & \text{for $ 0<x \le 1 $}, 
\\
I_{\nu}(x) \le C \e^x & \text{for $ x \ge 1 $}. 
\end{cases}
\nonumber
\end{align}

\begin{Lem} \label{rho}
For any $ u>0 $, $ a>0 $ and $ y>0 $, it holds that 
\begin{align}
\rho^a_u(y) \le \frac{2Cu}{a^2} 
 \quad \text{and}\quad
\tilde{\rho}^a_u(y) \le C . 
\label{eq: rho1}
\end{align}
For any fixed $ u>0 $ and $ y>0 $, it holds that 
\begin{align}
\tilde{\rho}^a_u(y) \tend{}{a \to \infty } 1 . 
\label{eq: rho2}
\end{align}
\end{Lem}

\Proof{
Using \eqref{eq: asympInu}, we easily have \eqref{eq: rho2}. \\
If $ 2 \sqrt{uy}/a \le 1 $, we have 
\begin{align}
\rho^a_u(y) \le C \frac{2u}{a^2} 
 \quad \text{and}\quad
\tilde{\rho}^a_u(y) \le C . 
\nonumber
\end{align}
If $ 2 \sqrt{uy}/a > 1 $, we have 
\begin{align}
\rho^a_u(y) 
\le& C \e^{-(\sqrt{u}+\sqrt{y})^2/a} \frac{\sqrt{u/y}}{a} 
\le C \frac{2 u}{a^2} , 
\nonumber \\
\tilde{\rho}^a_u(y) 
\le& C \e^{-(\sqrt{u}+\sqrt{y})^2/a} \le C . 
\nonumber
\end{align}
Therefore we obtain \eqref{eq: rho1}. 
}

\begin{Lem} \label{3--1}
Let $ f \in \mathcal{L}^1_+ $, $ x \in I $ and $ u>0 $. 
For any $ a \in I $, set 
\begin{align}
N_t^{a,u} =& a \bP_x \sbra{ f(L_{\eta^a_u}) 1_{\{ t<\eta^a_u \}} \mid \cF_t } , 
\nonumber \\
M_t^{a,u} =& a \bP_x \sbra{ f(L_{\eta^a_u}) \mid \cF_t } . 
\nonumber
\end{align}
Then it holds that 
$ N^{a,u}_t \to M^{s,f}_t $ and $ M^{a,u}_t \to M^{s,f}_t $ 
in probability with respect to $ \bP_x $ 
as $ a \to \infty $, 
where $ M^{s,f}_t $ has been defined in \eqref{eq: 2-7}. 
\end{Lem}

\Proof{
In what follows in this section we sometimes write to simplify $ M_t $ 
 for $ M^{s,f}_t $. 

(i) 
By the strong Markov property and by Lemma \ref{LawLeta2}, we have, for $ a>X_t $, 
\begin{align}
N^{a,u}_t 
= a \left. \bP_{X_t}[ f(b + L_{\eta^a_{u-c}}) ] 
\right|_{\begin{subarray}{l} b=L_t \\ c=L^a_t \end{subarray}} 
1_{\{ t<\eta^a_u \}} 
= {\rm (I)}_a + {\rm (II)}_a , 
\nonumber
\end{align}
where 
\begin{align}
{\rm (I)}_a 
=& X_t \left. \cbra{ \e^{- \frac{u-c}{a}} f(b) 
+ \int_0^{\infty } f(b + y) \rho^a_{u-c}(y) \d y 
} 
\right|_{\begin{subarray}{l} b=L_t \\ c=L^a_t \end{subarray}} 
1_{\{ t<\eta^a_u \}} , 
\nonumber \\
{\rm (II)}_a 
=& \left. \rbra{ 1 - \frac{X_t}{a} } \int_0^{\infty } f(b+y) \tilde{\rho}^a_{u-c}(y) \d y 
\right|_{\begin{subarray}{l} b=L_t \\ c=L^a_t \end{subarray}} 
1_{\{ t<\eta^a_u \}} . 
\nonumber
\end{align}
Letting $ a \to \infty $, we deduce from Lemma \ref{rho} that 
in probability with respect to $ \bP_x $ 
\begin{align}
{\rm (I)}_a 
\tend{}{a \to \infty }& X_t f(L_t) , 
\nonumber \\
{\rm (II)}_a 
\tend{}{a \to \infty }& \int_0^{\infty } f(L_t+y) \d y . 
\nonumber
\end{align}
We thus obtain 
$ N^{a,u}_t \to M^{s,f}_t $ in probability with respect to $ \bP_x $. 
Set 
\begin{align}
A^{a,u}_t = M^{a,u}_t - N^{a,u}_t 
= a f(L_{\eta^a_u}) 1_{\{ \eta^a_u \le t \}} . 
\nonumber
\end{align}
Since $ A^{a,u}_t \to 0 $, we obtain 
$ M^{a,u}_t \to M^{s,f}_t $ 
in probability with respect to $ \bP_x $. 
}

\begin{Thm}\label{theo:etaa}
Let $ f \in \mathcal{L}^1_+ $, $ x \in I $ and $ u>0 $. 
\begin{enumerate}
\item For any finite stopping time $ T $, there is the $\mathcal{L}^1$ convergence
\begin{align}
N^{a,u}_T \tend{}{a \to \infty } M^{s,f}_T 
\quad \text{in $ \mathcal{L}^1(\bP_x) $} . 
\nonumber
\end{align}
Consequently, 
for any bounded adapted process $ (F_t) $, it holds that 
\begin{align}
a \bP_x [F_T f(L_{\eta^a_u}) ; T<\eta^a_u] 
\tend{}{a \to \infty } \bP_x[F_T M^{s,f}_T] . 
\nonumber
\end{align}
\item Assume furthermore that $ \ell'$ $(=\ell=\infty) $ is 
either type-1-natural or type-2-natural. Then, for any bounded stopping time $ T $, we have
\begin{align}
M^{a,u}_T \tend{}{a \to \infty } M^{s,f}_T 
\quad \text{in $ \mathcal{L}^1(\bP_x) $} 
\label{eq: 3--4}
\end{align}
and, for any bounded adapted process $ (F_t) $, it holds that 
\begin{align}
a \bP_x \sbra{ F_T f(L_{\eta^a_u}) } 
\tend{}{a \to \infty } 
\bP_x \sbra{ F_T M^{s,f}_T } . 
\label{eq: 3--5}
\end{align}
\end{enumerate}
\end{Thm}

\Proof{(i) By the proof of Lemma \ref{3--1} 
and by Lemma \ref{rho}, we obtain, for $ a>1 $, 
\begin{align}
N^{a,u}_t 
\le& X_t f(L_t) 
+ \rbra{ \frac{2Cu}{a} + C } \int_0^{\infty } f(L_t+y) \d y 
\nonumber \\
\le& M^{s,f}_t 
+ \rbra{ 2Cu + C } \int_0^{\infty } f(y) \d y , 
\nonumber
\end{align}
where the last quantity is integrable with respect to $ \bP_x $. 
Thus we obtain the desired result by the dominated convergence theorem. 

(ii) Observe that since $ (M^{s,f}_t) $ is a $ \bP_x $-martingale, 
we may apply the optional stopping theorem to obtain  
\begin{align}
\bP_x[A^{a,u}_T] 
=& \bP_x \sbra{ X_{\eta^a_u} f(L_{\eta^a_u}) ; \eta^a_u \le T } 
\nonumber \\
\le& \bP_x \sbra{ M_{\eta^a_u} ; \eta^a_u \le T } 
\nonumber \\
=& \bP_x \sbra{ M_T ; \eta^a_u \le T } 
\tend{}{a \to \infty } 0 . 
\nonumber
\end{align}
Since $ A^{a,u}_T \to 0 $, $ \bP_x $-a.s., 
we see that $ A^{a,u}_T \to 0 $ in $ \mathcal{L}^1(\bP_x) $. 
Hence we obtain \eqref{eq: 3--4} and \eqref{eq: 3--5}. 
}



\subsection{Limit as $ u $ tends to infinity with $ a $ being fixed}

Suppose $ \ell'$ $(=\ell=\infty) $ is 
either entrance, type-1-natural or type-2-natural. 
We thus have, for any $ x,a \in I $, 
\begin{align}
\bP_x(\eta^a_u<\infty ) = 1 
\quad \text{and} \quad 
\eta^a_u \tend{}{u \to \infty } \infty 
\ \text{$ \bP_x $-a.s.} 
\nonumber
\end{align}
In fact, $ \eta^a_u $ increases to a limit $ \eta^a_{\infty } $ 
which must be infinite $ \bP_x $-a.s. by \eqref{eq: Eqeta}. 
For the clock $ (\tau = \eta^a_u, \,u>0 )$, 
we only consider the weights $ f(L_{\eta^a_u}) $ 
for $ f(u) = \e^{- \beta u} $ and $ f(u) = 1_{\{ u=0 \}} $.

\begin{Lem} \label{ulim1}
Let $ x,a \in I $, $ \beta>0 $ and $ t>0 $. 
For $ u>0 $, set 
\begin{align}
N^{u,\beta,a}_t = \e^{ \frac{\beta u}{1+\beta a} } 
\bP_x \sbra{ \left. \e^{- \beta L_{\eta^a_u}} 1_{\{ t<\eta^a_u \}} \right| \cF_t } , 
\quad 
M^{u,\beta,a}_t = \e^{ \frac{\beta u}{1+\beta a} } 
\bP_x \sbra{ \left. \e^{- \beta L_{\eta^a_u}} \right| \cF_t } 
\nonumber
\end{align}
and 
\begin{align}
M^{\beta,a}_t =& \frac{1 + \beta (X_t \wedge a)}{1+\beta a} 
\exp \rbra{- \beta L_t + \frac{\beta}{1+\beta a} L^a_t} . 
\nonumber
\end{align}
Then it holds that 
$ N^{u,\beta,a}_t \to M^{\beta,a}_t $ and $ M^{u,\beta,a}_t \to M^{\beta,a}_t $, $ \bP_x $-a.s. 
as $ u \to \infty $. 
\end{Lem}

\Proof{
By the strong Markov property and by Lemmas \ref{LawLeta} and \ref{LawLeta2}, 
we have, for $ u $ large enough to have $ \eta^a_u > t $, 
\begin{align}
N^{u,\beta,a}_t 
=& \e^{ \frac{\beta u}{1+\beta a} } \exp \rbra{- \beta L_t} 
\left. \bP_{X_t} \sbra{ \exp \rbra{ - \beta L_{\eta^a_{u-c}}) } } 
\right|_{c=L^a_t} 
\nonumber \\
=& 
\cbra{ 1_{\{ a \le X_t \}} + \frac{1+\beta X_t}{1+\beta a} 1_{\{ X_t<a \}} } 
\exp \rbra{ - \beta L_t + \frac{\beta}{1+\beta a} L^a_t } 
\nonumber \\
=& 
M^{\beta,a}_t . 
\nonumber
\end{align}
Thus we obtain 
$ N^{u,\beta,a}_t \to M^{\beta,a}_t $, $ \bP_x $-a.s. 
as $ u \to \infty $. 
Since 
\begin{align}
A^{u,\beta,a}_t 
:= M^{u,\beta,a}_t - N^{u,\beta,a}_t 
= \e^{ \frac{\beta u}{1+\beta a} } \e^{- \beta L_{\eta^a_u}} 1_{\{ \eta^a_u \le t \}} , 
\nonumber
\end{align}
we have $ A^{u,\beta,a}_t \to 0 $, $ \bP_x $-a.s., 
and thus we obtain 
$ M^{u,\beta,a}_t \to M^{\beta,a}_t $, $ \bP_x $-a.s. 
}

\begin{Thm} \label{ulim2}
Let $ x,a \in I $ and $ \beta>0 $. 
Then, for any $ t>0 $, it holds that 
\begin{align}
N^{u,\beta,a}_t \tend{}{u \to \infty } M^{\beta,a}_t 
\quad \text{and} \quad 
M^{u,\beta,a}_t \tend{}{u \to \infty } M^{\beta,a}_t 
\quad 
\text{in $ \mathcal{L}^1(\bP_x) $}. 
\nonumber
\end{align}
Consequently, for any bounded adapted process $ (F_t) $, 
it holds that 
\begin{align}
\lim_{u \to \infty } 
\e^{ \frac{\beta u}{1+\beta a} } 
\bP_x[F_t \e^{-\beta L_{\eta^a_u}} ; t < \eta^a_u] 
= \lim_{u \to \infty } 
\e^{ \frac{\beta u}{1+\beta a} } 
\bP_x[F_t \e^{-\beta L_{\eta^a_u}}] 
= \bP_x[F_t M^{\beta,a}_t] . 
\nonumber
\end{align}
It also holds that $ (M^{\beta,a}_t) $ is a $ \bP_x $-martingale. 
\end{Thm}

\Proof{
Let us first prove that 
$ \bP_x[\e^{c L^a_t}] < \infty $ for all $ c>0 $ and $ t>0 $. 
Following the same argument as in the proof of Lemma \ref{lem1}, we obtain 
\begin{align}
\bP_a \sbra{ \exp \rbra{ c L^a_{\ve_q} } } 
= \frac{1}{r_q(a,a)} \int_0^{\infty } \e^{c u} \e^{-u/r_q(a,a)} \d u . 
\nonumber
\end{align}
Since $ r_q(a,a) \to 0 $ as $ q \to \infty $, 
we may take $ q>0 $ large enough so that $ r_q(a,a) < 1/c $. 
This shows that $ \bP_a \sbra{ \exp \rbra{ c L^a_{\ve_q} } } < \infty $. 
By the monotonicity, we see that 
$ \bP_x[\e^{c L^a_t}] < \infty $ for all $ t>0 $. 

The fact that $ L^a_t $ admits exponential moments implies 
that $ M^{\beta,a}_t \in \mathcal{L}^1(\bP_x) $ for all $ t>0 $. 
Thus, by the dominated convergence theorem, 
we see that $ N^{u,\beta,a}_t \tend{}{u \to \infty } M^{\beta,a}_t $ 
in $ \mathcal{L}^1(\bP_x) $ for all $ t>0 $. 

\noindent
We second note that, for $ q>0 $, 
\begin{align}
\bP_x(\eta^a_u \le t) 
\le \e^{qt} \bP_x[\e^{-q \eta^a_u}] 
\le \e^{qt} \bP_a[\e^{-q \eta^a_u}] 
= \e^{qt} \e^{- u/r_q(a,a)} . 
\nonumber
\end{align}
We may take $ q>0 $ large enough so that $ r_q(a,a) < (1+\beta a)/\beta $. 
Then we obtain 
\begin{align}
\bP_x[A^{u,\beta,a}_t] 
\le \e^{\frac{\beta u}{1+\beta a}} \bP_x(\eta^a_u \le t) 
\le \e^{qt} \exp \cbra{ - \rbra{ \frac{1}{r_q(a,a)} - \frac{\beta}{1+\beta a} } u } 
\tend{}{u \to \infty } 0. 
\nonumber
\end{align}
Thus we obtain 
$ A^{u,\beta,a}_t \tend{}{u \to \infty } 0 $ 
in $ \mathcal{L}^1(\bP_x) $ for all $ t>0 $, 
which implies 
$ M^{u,\beta,a}_t \tend{}{u \to \infty } M^{\beta,a}_t $ 
in $ \mathcal{L}^1(\bP_x) $ for all $ t>0 $. 
}

\noindent
We conclude this section by looking at the weight $1_{\{L_{\eta_u^a =0}\}} = 1_{\{\eta_u^a < T_0\}}$.

\begin{Thm}
Let $ x,a \in I $. 
For $ u>0 $ and $ t>0 $, set 
\begin{align}
N^{u,\infty ,a}_t =& \e^{ u/a } 
\bP_x (t<\eta^a_u<T_0 \mid \cF_t) 
\nonumber \\
M^{u,\infty ,a}_t =& \e^{ u/a } 
\bP_x (\eta^a_u<T_0 \mid \cF_t) 
\nonumber \\
M^{\infty ,a}_t =& \frac{X_t \wedge a}{a} 
\e^{L^a_t/a} 1_{\{ t<T_0 \}} . 
\nonumber
\end{align}
Then,
\begin{align}
N^{u,\infty ,a}_t \tend{}{u \to \infty } M^{\infty ,a}_t 
\quad \text{and} \quad 
M^{u,\infty ,a}_t \tend{}{u \to \infty } M^{\infty ,a}_t 
\quad 
\text{$ \bP_x $-a.s. and in $ \mathcal{L}^1(\bP_x) $}. 
\nonumber
\end{align}
Consequently, for any bounded adapted process $ (F_t) $, 
it holds that 
\begin{align}
\lim_{u \to \infty } 
\e^{ u/a } 
\bP_x[F_t ; t < \eta^a_u < T_0] 
= \lim_{u \to \infty } 
\e^{ u/a } 
\bP_x[F_t ; \eta^a_u < T_0] 
= \bP_x[F_t M^{\infty ,a}_t] . 
\nonumber
\end{align}
It also holds that $ (M^{\infty ,a}_t) $ is a $ \bP_x $-martingale. 
\end{Thm}

\Proof{
Letting $ \beta \to \infty $, we see, 
from Lemma \ref{ulim1} and Theorem \ref{ulim2}, that 
\begin{align}
N^{u,\infty ,a}_t = M^{\infty ,a}_t 1_{\{ t< \eta^a_u \}} 
\nonumber
\end{align}
and 
\begin{align}
A^{u,\infty ,a}_t 
:= M^{u,\infty ,a}_t - N^{u,\infty ,a}_t 
= \e^{u/a} 1_{\{ \eta^a_u < T_0 \}} 1_{\{ \eta^a_u \le t \}} . 
\nonumber
\end{align}
The remainder of the proof 
is the same as that of Theorem \ref{ulim2}. 
}

\section{Universal $ \sigma $-finite measures} \label{uni}

In this section we shall describe the law of some penalized processes using universal $ \sigma $-finite measures. 
We deal with the transient and recurrent cases separately.

\subsection{The transient case}

\begin{Thm} \label{3-0}
Suppose $ \ell<\infty $, i.e., 0 is transient. 
Let $ f \in \mathcal{L}^1_+ $ and $ x \in I $. 
Let $ t $ be a constant time and let $ F_t $ be a bounded $ \cF_t $-measurable functional. 
Then 
\begin{align}
\lim_{q \down 0} \bP_x[F_t f(L_{\ve_q}) ; t<\ve_q] 
= \lim_{q \down 0} \bP_x[F_t f(L_{\ve_q})] 
= \bP_x[F_t f(L_{\infty })] . 
\label{eq: 3-0-1}
\end{align}
If, in particular, $ \ell $ is type-3-natural (see Section 7), then 
\begin{align}
\lim_{a \up \ell} \bP_x[F_t f(L_{T_a}) ; t<T_a ] 
= \lim_{a \up \ell} \bP_x[F_t f(L_{T_a})] 
= \bP_x[F_t f(L_{\infty })] . 
\label{eq: 3-0-2}
\end{align}
\end{Thm}

\Proof{
By Theorem \ref{2}, we see that 
\eqref{eq: 3-0-1} is equivalent to 
\begin{align}
\bP_x[F_t f(L_{\infty })] 
= \bP_x[F_t M_t] , 
\label{eq: 3-0-3}
\end{align}
where 
\begin{align}
M_t = \frac{1}{\ell} \cbra{ X_t f(L_t) 
+ \rbra{ 1-\frac{X_t}{\ell} } \int_0^{\infty } \e^{-u/\ell} f(L_t+u) \d u } . 
\nonumber
\end{align}
On the other hand, we use (i) of Theorem \ref{1} and obtain 
\begin{align}
\bP_x[f(L_{\infty })|\cF_t] 
= \left. \bP_{X_t}[f(a+L_{\infty })] \right|_{a=L_t} 
= M_t . 
\nonumber
\end{align}
Thus we obtain \eqref{eq: 3-0-3}. 

Using Theorem \ref{2-1} 
instead of Theorem \ref{2}, 
we can obtain \eqref{eq: 3-0-2} in the same way as above. 
}

\begin{Rem}
Observe that both penalizations yield the same measure $\bQ_x$ 
\begin{align}
\bQ_x = \frac{ f(L_\infty)}{\bP_x[f(L_\infty )]} \centerdot \bP_x, 
\nonumber
\end{align}
which is absolutely continuous with respect to $\bP_x$. This is not very surprising since the initial process spends little time at 0, hence the penalization by the local time at 0 has no real impact on the measure $\bP_x$.
\end{Rem}

\subsection{The recurrent case}

Let $ \bP^{(u)}_{x,y} $ denote the law of the bridge 
with duration $ u $ starting from $ x $ and ending at $ y $. 
Following \cite{MR1278079}, this measure can be characterized by 
\begin{align}
\bP^{(u)}_{x,y}(A) = \bP_x \sbra{ 1_A \frac{p_{u-t}(X_t,y)}{p_u(x,y)} } 
, \quad A \in \cF_t , \ 0<t<u , 
\nonumber
\end{align}
where $ p_u(x,y) $ denotes 
the transition density of the process $ X $ with respect to $ \d m(y) $. 
We have the conditioning formula: 
\begin{align}
\bP_x \sbra{ \int_0^{\infty } F_u \d L_u } 
= \int_0^{\infty } \bP_x[\d L_u] \bP^{(u)}_{x,0}[F_u] 
\nonumber
\end{align}
for all non-negative predictable processes $ (F_u) $, 
where we write symbolically (see It\^o--McKean \cite[p.183]{IMK}):
\begin{align}
\bP_x[\d L_u] = p_u(x,0) \d u . 
\nonumber
\end{align}
We also have the last exit decomposition formula (see Biane-Yor \cite{BY}): 
\begin{align}
\bP_x[F_t;T_0 \le t] 
= \int_0^t \bP_x[\d L_u] \rbra{ \bP^{(u)}_{x,0} \bullet \vn^{[t-u]} }[F_t] 
\nonumber
\end{align}
for all non-negative $ \cF_t $-measurable functionals $ (F_t) $, 
where $\bullet$ denotes the concatenation operator and  
\begin{align}
\vn^{[t]}(\cdot) = \vn(\cdot \cap \{ t<T_0 \}) . 
\nonumber
\end{align}

\noindent
For $ h=h_0 $ (see \eqref{eq:h0}) or $ h=s $ (the scale function), let $ \bP^h_x $ denote the law of the $ h $-transform: 
\begin{align}
\bP^h_x(A;t<\zeta) =& \frac{1}{h(x)} \bP_x[1_{A\cap\{t<T_0\}} h(X_t)] 
\quad (x>0) , 
\nonumber \\
\bP^h_0(A;t<\zeta) =& \vn[1_A h(X_t)] 
\nonumber
\end{align}
for $ A \in \cF_t $ and where $\zeta$ denotes the lifetime of the process. 
Note that, when $ h=h_0 $ or $ h=s $, the coordinate process under $ \bP_x^h $ 
never hits zero; see \cite[Theorems 7.6 and 7.3]{YYrenorm}. 
We now define the $\sigma$-finite measure 
\begin{align}
\cP^h_x = \int_0^{\infty } \bP_x[\d L_u] \rbra{ \bP^{(u)}_{x,0} \bullet \bP^h_0 } 
+ h(x) \bP^h_x . 
\label{eq:defcP}
\end{align}

\begin{Thm} \label{3-1}
Suppose $ \ell=\infty $, i.e., 0 is recurrent. 
Let $ f \in \mathcal{L}^1_+ $ and $ x \in I $. 
Let $ t $ be a constant time and let $ F_t $ be a bounded $ \cF_t $-measurable functional. 
Then 
\begin{align}
\lim_{q \down 0} H(q) \bP_x[F_t f(L_{\ve_q}) ; t<\ve_q] 
= \cP^{h_0}_x[F_t f(L_{\zeta});t<\zeta] . 
\nonumber
\end{align}
\end{Thm}

\Proof{
By Theorem \ref{2}, it suffices to show 
\begin{align}
\cP^{h_0}_x[F_t f(L_{\zeta});t<\zeta] 
= \bP_x\left[F_t N^{h_0,f}_t\right] . 
\label{eq: 3-1-1}
\end{align}

\noindent
Denote $ g = \sup \{ t < \zeta: X_t=0 \} $, where $ \sup \emptyset = 0 $. 
Observe first that on the set $ \{ 0 = g \le t<\zeta \} $, we have 
\begin{align}
\cP^{h_0}_x[F_t f(L_{\zeta});0 = g \le t<\zeta] 
=& h_0(x) \bP^{h_0}_x[F_t f(L_t) ; t<\zeta] 
\nonumber \\
=& \bP_x[F_t f(L_t) h_0(X_t); t<T_0] . 
\nonumber
\end{align}
Next, on the set $ \{ 0 < g \le t<\zeta \} $, we have 
\begin{align}
\cP^{h_0}_x[F_t f(L_{\zeta});0 < g \le t<\zeta] 
=& \int_0^t \bP_x[\d L_u] \rbra{ \bP^{(u)}_{x,0} \bullet \bP^h_0 }[F_t f(L_t);t<\zeta] 
\nonumber \\
=& \int_0^t \bP_x[\d L_u] \rbra{ \bP^{(u)}_{x,0} \bullet \vn }[F_t f(L_t) h_0(X_t)] 
\nonumber \\
=& \bP_x[F_t f(L_t) h_0(X_t) ; T_0 \le t] . 
\nonumber
\end{align}
Finally, on the set $ \{ t < g < \zeta \} $, we have 
\begin{align}
\cP^{h_0}_x[F_t f(L_{\zeta}); t < g<\zeta] 
=& \int_t^{\infty } \bP_x[\d L_u] \bP^{(u)}_{x,0}[F_t f(L_u)] 
\nonumber \\
=& \bP_x \sbra{ F_t \int_t^{\infty } f(L_u) \d L_u } 
\nonumber \\
=& \bP_x \sbra{ F_t \int_{L_t}^{\infty } f(u) \d u } . 
\nonumber
\end{align}
Summing all three terms yields \eqref{eq: 3-1-1}. 
}

\begin{Thm} \label{3-2}
Suppose $ \ell' $ is either entrance, type-1-natural or type-2-natural. 
Let $ f \in \mathcal{L}^1_+ $ and $ x \in I $. 
Let $ t $ be a constant time and let $ F_t $ be a bounded $ \cF_t $-measurable functional. 
Then 
\begin{align}
\lim_{a \up \ell} a \bP_x[F_t f(L_{T_a}) ; t<T_a ] 
= \cP^s_x[F_t f(L_{\zeta}) ; t<\zeta] . 
\nonumber
\end{align}
\end{Thm}

The proof is similar to that of Theorem \ref{3-1}, 
where we use Theorem \ref{2-1} instead of Theorem \ref{2}, so we omit it. \\

\begin{Rem}
Define, for $h=h_0$ or $h=s$, the penalized measures
$$\bQ_{x}^{h,f}(A ; t<\zeta) =  \bP_{x}\left[1_{A} \frac{N_t^{h,f}}{N_0^{h,f}}\right]\quad\qquad (\text{for }A\in \cF_t).$$
Looking at (\ref{eq:defcP}), we see that, under the assumptions of Theorems \ref{3-1} or \ref{3-2}, the paths of the coordinates processes under $\bQ_{x}^{h,f}$ are essentially given, up to some killing time $\zeta$, by the concatenation of a weighted bridge of the original diffusion, and a process conditioned not to hit 0. In particular, the penalized process is no longer recurrent, even if $\zeta=\infty$, i.e. if $N^{h,f}$ is a $\bP_x$-martingale. Of course, in this case, the two probability measures $\bQ_{x}^{h,f}$ and $\bP_x$ are singular.

%
\end{Rem}

\section{Exponential weights} \label{ex}

Let us investigate the example where we take 
\begin{align}
f(x) = \e^{- cx} 
, \quad c>0 . 
\nonumber
\end{align}
In this specific case, the penalized process remains a generalized diffusion, which is not the case with other functions $f$. 
The supermartingales $ N_t = N^{h_0,f}_t $ and $ N_t = M^{s,f}_t $, 
which have been given by \eqref{eq: Nh0ft} and \eqref{eq: 2-7}, respectively, 
may be represented at the same time as 
\begin{align}
N_t = h^c(X_t) \e^{-cL_t} 
\nonumber
\end{align}
where $ h = h_0 $ and $ s $, respectively, and 
\begin{align}
h^c(x) = h(x) + \frac{1 - \frac{x}{\ell}}{c + \frac{1}{\ell}} . 
\nonumber
\end{align}
Since $ (N_t) $ is a supermartingale, 
we may define the subprobability measure $ \bQ^{h,c}_x $ by 
\begin{align}
\bQ^{h,c}_x(A;t<\zeta) 
= \bP_x \sbra{ \frac{h^c(X_t)}{h^c(x)} \e^{-cL_t} ;A } 
\quad \text{for $ A \in \cF_t $ and $ t \ge 0 $} . 
\nonumber
\end{align}
Then the process $ \{ X,(\bQ^{h,c}_x)_{x \in I} \} $ is a diffusion on $ I $ 
whose local generator on $ (0,\ell') $ without killing part 
is given as $ (h^c)^{-1} \frac{\d}{\d m} \frac{\d}{\d s} h^c $. 
Thus the corresponding speed measure and  scale function are given as 
\begin{align}
m^{h,c}(x) = \int_{(0,x]} h^c(y)^2 \d m(y) 
, \quad 
s^{h,c}(x) = \int_0^x \frac{\d y}{h^c(y)^2} . 
\nonumber
\end{align}
Denote  $\displaystyle \rho_q = \phi_q - \frac{\psi_q}{H(q)}$ and
\begin{align}
\phi_q^{h,c} = h^c(0) \cdot \frac{\phi_q + c \psi_q}{h^c} 
, \quad 
\rho_q^{h,c} = h^c(0) \cdot \frac{\rho_q}{h^c} . 
\nonumber
\end{align}
Then we obtain (see Theorems 7.3 and 7.6 of \cite{YYrenorm}) that $ \varphi = \phi_q^{h,c} $ (resp. $ \rho_q^{h,c} $) 
is a positive increasing (resp. decreasing) solution to the differential equation 
\begin{align}
\rbra{ \frac{\d }{\d {m^{h,c}}} \frac{\d }{\d {s^{h,c}}} - \frac{\pi_0}{h^c} } \varphi = q \varphi 
\ \text{(if $ h=h_0 $)}, 
\quad 
\rbra{ \frac{\d }{\d {m^{h,c}}} \frac{\d }{\d {s^{h,c}}} } \varphi = q \varphi 
\ \text{(if $ h=s $)} 
\nonumber
\end{align}
which satisfies the boundary condition 
\begin{align}
\phi_q^{h,c}(0) = 1 
\quad \text{and} \quad 
\frac{\d \phi_q^{h,c}}{\d s^{h,c}}(0) = 0 
\quad \text{(resp. $ \rho_q^{h,c}(0) = 1 $)}. 
\nonumber
\end{align}
Note that we have used here the values 
\begin{align}
h^c(0) = \frac{1}{c+\frac{1}{\ell}} 
 \quad \text{and}\quad 
(h^c)'(0) = \frac{c}{c+\frac{1}{\ell}} . 
\nonumber
\end{align}

\begin{Thm}
The resolvent operator for the diffusion $ \{ X,(\bQ^{h,c}_x)_{x \in I} \} $ is given as 
\begin{align}
\bQ^{h,c}_x \sbra{\int_0^{\infty } \e^{-qt} f(X_t) \d t } 
= \int_I r^{h,c}_q(x,y) f(y) \d m^{h,c}(y) 
, \quad q>0 , 
\nonumber
\end{align}
where 
\begin{align}
r^{h,c}_q(x,y) = r^{h,c}_q(y,x) 
= \frac{H(q)}{h^c(0)^2(cH(q) + 1)} \phi^{h,c}_q(x) \rho^{h,c}_q(y) 
, \quad x,y \in I , \ x \le y . 
\label{eq: 4-1}
\end{align}
Consequently, 0 for $ \{ X,(\bQ^{h,c}_x)_{x \in I} \} $ is regular-reflecting. 
\end{Thm}

\Proof{
Let $ \varphi^c(x) = \varphi(x) h^c(x) $. Then we have 
\begin{align}
& \bP_x \sbra{ \int_0^{\infty } \e^{-qt} \varphi^c(X_t) \e^{-cL_t} \d t } 
\nonumber \\
=& \bP_x \sbra{ \int_0^{T_0} \e^{-qt} \varphi^c(X_t) \d t } 
+ \bP_x[\e^{-qT_0}] \bP_0 \sbra{ \int_0^{\infty } \e^{-qt} \varphi^c(X_t) \e^{-cL_t} \d t } 
\nonumber \\
=& R^0_q \varphi^c(x) 
+ \bP_x[\e^{-qT_0}] 
\bP_0 \sbra{ \sum_u \e^{-cu-q\eta_{u-}^0} \int_0^{T_0(p(u))} \e^{-qt} \varphi^c(p(u)_t) \d t } 
\nonumber \\
=& R^0_q \varphi^c(x) 
+ \bP_x[\e^{-qT_0}] \bP_0 \sbra{ \int_0^{\infty } \e^{-cu-q\eta_{u}^0} \d u } 
\vn \sbra{ \int_0^{T_0} \e^{-qt} \varphi^c(X_t) \d t } 
\nonumber \\
=& R^0_q \varphi^c(x) 
+ \bP_x[\e^{-qT_0}] \cdot \frac{1}{c+\frac{1}{H(q)}} \cdot \frac{R_q \varphi^c(0)}{H(q)} . 
\nonumber
\end{align}
Since $ \bP_x[\e^{-qT_0}] R_q \varphi^c(0) = R_q \varphi^c(x) - R^0_q \varphi^c(x) $, we obtain 
\begin{align}
\bQ^{h,c}_x \sbra{ \int_0^{\infty } \e^{-qt} \varphi(X_t) \d t } 
= \frac{1}{h^c(x)} \cbra{ \frac{1}{cH(q)+1} R_q\varphi^c(x) 
+ \frac{cH(q)}{cH(q)+1} R^0_q\varphi^c(x) } . 
\nonumber
\end{align}
From this we obtain \eqref{eq: 4-1}. 
}

\begin{Rem}
The boundary classification at $ \ell' $ 
is the same as that for the $ h $-transform of the stopped process; 
see Theorems 7.3 and 7.6 of \cite{YYrenorm}. 
\end{Rem}

\section{Appendix: the boundary classification} \label{app}

The following tables explain the boundary classification which we take from \cite{YYrenorm} 
and the recurrence property of the corresponding diffusion to each class: 

\begin{center}
\begin{tabular}{lllll}
$ x=\ell' $ & & $ I' $ & $ I $ & $ x=0 $ 
\\ \hline
regular-reflecting 	& $ \ell'<\ell=\infty $ 	& $ [0,\ell'] $ & $ =I' $ & positive recurrent 
\\
regular-elastic 	& $ \ell'<\ell<\infty $ 	& $ [0,\ell'] $ & $ [0,\ell'] \cup \{ \ell \} $ & transient 
\\
regular-absorbing 	& $ \ell'=\ell<\infty $ 	& $ [0,\ell) $ & $ [0,\ell] $ & transient 
\\
exit 				& $ \ell'=\ell<\infty $ 	& $ [0,\ell) $ & $ [0,\ell] $ & transient 
\\
entrance 			& $ \ell'=\ell=\infty $ 	& $ [0,\infty ) $ & $ =I' $ & positive recurrent 
\\
type-1-natural 		& $ \ell'=\ell=\infty $	& $ [0,\infty ) $ & $ =I' $ & null recurrent 
\\
type-2-natural 		& $ \ell'=\ell=\infty $	& $ [0,\infty ) $ & $ =I' $ & positive recurrent 
\\
type-3-natural 		& $ \ell'=\ell<\infty $ 	& $ [0,\ell) $ & $ =I' $ & transient 
\end{tabular}
\end{center}

\begin{center}
\begin{tabular}{l||l|l}
& $ \ell=\infty $ & $ \ell<\infty $ 
\\ \hline\hline
$ m(\infty )=\infty $ & {\bf (1)} 0 is null-recurrent & {\bf (3)} 0 is transient 
\\
$ \pi_0=0 $ & $ [\ell'=\ell=\infty] $ & $ [\ell'<\ell<\infty] $ 
\\
& $ \ell' $ is type-1-natural & $ \ell' $ is regular-elastic 
\\
& & $ [\ell'=\ell<\infty] $ 
\\
& & $ \ell' $ is regular-absorbing 
\\
& & \hspace{6.6mm} exit 
\\
& & \hspace{6.6mm} type-3-natural 
\\ \hline
$ m(\infty )<\infty $ & {\bf (2)} 0 is positive recurrent & {[impossible]} 
\\
$ \pi_0>0 $ & $ [\ell'<\ell=\infty ] $ 
\\
& $ \ell' $ is regular-reflecting 
\\
& $ [\ell'=\ell=\infty ] $ 
\\
& $ \ell' $ is entrance 
\\
& \hspace{6.6mm} type-2-natural 
\end{tabular}
\end{center}

As the reader may not be familiar with our classification of boundaries, 
it may be useful to give below some examples of computation of boundaries. 
Let $ \tilde{X} $ be a diffusion on $ [0,\infty ) $ where 0 is the reflecting boundary 
and whose local generator on $ (0,\infty ) $ is given by 
\begin{align}
\tilde{L}f = \frac{1}{2}(f'' - b f') = \frac{\d}{\d \tilde{m}} \frac{\d}{\d \tilde{s}} f 
\quad \text{on $ C_c((0,\infty )) $} 
\nonumber
\end{align}
for some function $ b $ of the form $ b(x) = cx^{\nu-1} $, 
which we may call the {\em power drift}. 
Then its scale change $ X = \tilde{s}(\tilde{X}) $ is a diffusion with natural scale $ s(x)=x $ 
and with speed measure $ \d m(x) $ defined by $ m = \tilde{m} \circ \tilde{s}^{-1} $. 
\begin{enumerate}
\item 
Let $ \alpha $ be a constant and 
\begin{align}
\tilde{L}f = \frac{1}{2}f'' - \frac{2 \alpha -1}{2x} f' 
= \frac{\d}{\d \tilde{m}} \frac{\d}{\d \tilde{s}} f 
\quad \text{on $ C_c((0,\infty )) $} , 
\nonumber
\end{align}
where we may choose $ \tilde{m}(x) = \frac{2}{2-2\alpha } x^{2-2\alpha } $ 
and $ \tilde{s}(x) = \frac{1}{2 \alpha } x^{2 \alpha } $. 
The corresponding diffusion is called the {\em reflecting Bessel process of index $ \alpha $}. 
As we require that $ 0 $ is regular-reflecting, 
we assume $ 0 < \alpha < 1 $. 
If we take $ m = \tilde{m} \circ \tilde{s}^{-1} $, 
then it falls into the case {\bf (1)} above. 

\item 
Let $ c $ and $ \nu $ be non-zero constants and 
\begin{align}
\tilde{L} f = \frac{1}{2} \rbra{ f'' - c \nu x^{\nu-1} f' } 
\quad \text{on $ C_c((0,\infty )) $}. 
\nonumber
\end{align}
If $ \nu=1 $, then it is a Brownian motion with constant negative drift. 
If $ \nu=2 $, then it is an Ornstein--Uhlenbeck process. 
As we require that $ 0 $ is regular-reflecting, 
we assume $ c>0 $ and $ \nu>0 $. 
In this case we may choose 
\begin{align}
s' = \e^{c x^{\nu}} 
, \quad 
s = \int_0^x \e^{c y^{\nu}} \d y 
\nonumber
\end{align}
and 
\begin{align}
m' = 2 \e^{- c x^{\nu}} 
, \quad 
m = 2 \int_0^x \e^{- c y^{\nu}} \d y . 
\nonumber
\end{align}
In particular, we have $ m(\infty ) < \infty $. 
Note that 
\begin{align}
J:= \frac{1}{2} \int_1^{\infty } \{ s(x)-s(1) \} \d m(x) 
=& \int_1^{\infty } \rbra{ \int_1^x \e^{cy^{\nu}} \d y } \e^{-cx^{\nu}} \d x 
\nonumber \\
=& \int_1^{\infty } \rbra{ \int_y^{\infty } \e^{-cx^{\nu}} \d x } \e^{cy^{\nu}} \d y . 
\nonumber
\end{align}
We shall prove that 
\begin{align}
\text{$ \infty $ is } 
\begin{cases}
\text{type-2-natural if $ 0 < \nu \le 2 $}, \\
\text{entrance if $ 2 < \nu < \infty $}
\end{cases}
\nonumber
\end{align}
which is equivalent as saying that $J=+\infty$ (resp. $J<+\infty$), see \cite{YYrenorm}.

If $ 1 \le \nu \le 2 $, then 
\begin{align}
\int_1^x \e^{cy^{\nu}} \d y 
=& \int_1^x (\e^{cy^{\nu}})' \frac{y^{1-\nu}}{c \nu} \d y 
\nonumber \\
=& \sbra{ \e^{cy^{\nu}} \frac{y^{1-\nu}}{c \nu} }_1^x 
+ \frac{\nu-1}{c \nu} \int_1^x \e^{cy^{\nu}} y^{-\nu} \d y 
\nonumber \\
\ge& \e^{cx^{\nu}} \frac{x^{1-\nu}}{c \nu} - c' 
\nonumber
\end{align}
for some constant $ c'>0 $. Hence we have 
\begin{align}
J \ge \frac{1}{c \nu} \int_1^{\infty } x^{1-\nu} \d x 
- c' \int_1^{\infty } \e^{- c x^{\nu}} \d x = \infty . 
\nonumber
\end{align}

For $ \nu>0 $, we have 
\begin{align}
\int_y^{\infty } \e^{-cx^{\nu}} \d x 
=& - \int_y^{\infty } (\e^{-cx^{\nu}})' \frac{x^{1-\nu}}{c \nu} \d y 
\nonumber \\
=& - \sbra{ \e^{-cx^{\nu}} \frac{x^{1-\nu}}{c \nu} }_y^{\infty } 
+ \frac{1-\nu}{c \nu} \int_y^{\infty } \e^{-cx^{\nu}} x^{-\nu} \d x . 
\nonumber
\end{align}
If $ 0 < \nu < 1 $, then 
\begin{align}
\int_y^{\infty } \e^{-cx^{\nu}} \d x 
\ge& \e^{-cy^{\nu}} \frac{y^{1-\nu}}{c \nu} 
\quad \text{and} \quad 
J \ge \frac{1}{c\nu} \int_1^{\infty } y^{1-\nu} \d y = \infty . 
\nonumber
\end{align}
If $ \nu > 2 $, then 
\begin{align}
\int_y^{\infty } \e^{-cx^{\nu}} \d x 
\le \e^{-cy^{\nu}} \frac{y^{1-\nu}}{c \nu} 
\quad \text{and} \quad 
J \le \frac{1}{c\nu} \int_1^{\infty } y^{1-\nu} \d y < \infty . 
\nonumber
\end{align}
\end{enumerate}

\def\cprime{$'$} \def\cprime{$'$}

\end{document}